\documentclass{article}
\usepackage{amssymb}
\usepackage{amsfonts}
\usepackage{amsmath,amsthm,amssymb,latexsym,amscd}
\usepackage{epsf}
\usepackage[pdftex]{graphicx}
\newtheorem{lemma}{Lemma}
\newtheorem{theorem}{Theorem}
\newtheorem{proposition}{Proposition}
\newtheorem{remark}{Remark}
\newtheorem{definition}{Definition}
\newtheorem{hypothesis}{Hypothesis}
\title{Dynamic Brittle Fracture as a  Small Horizon Limit of Peridynamics}
\author{Robert Lipton\thanks{Department of Mathematics,
        Louisiana State University,
        Baton Rouge, LA 70803,
        ({\tt lipton@math.lsu.edu}).}}
\begin{document}
\maketitle
\begin{abstract} 
We consider the nonlocal formulation of continuum mechanics described by peridynamics. 
We provide a link between peridynamic evolution and brittle fracture evolution for a broad class of peridynamic potentials associated with unstable peridynamic constitutive laws. Distinguished limits of peridynamic evolutions are identified that correspond  to vanishing peridynamic horizon. The limit evolution is associated with dynamic brittle fracture and satisfies a dynamic energy  inequality expressed in terms of the kinetic energy of the motion together with a bulk elastic energy and a Griffith surface energy. It corresponds to the simultaneous evolution of elastic displacement and brittle fracture with displacement fields satisfying the wave equation inside the cracking domain. The wave equation provides the dynamic coupling between elastic waves and the evolving fracture path inside the media. The elastic moduli, wave speed and energy release rate for the evolution are explicitly determined by moments of the peridynamic influence function and the peridynamic potential energy. 
\end{abstract}
\begin{flushleft} 
\end{flushleft}
%
%
%
\pagestyle{myheadings}
\markboth{R. LIPTON}{Peridynamics and the Small Horizon Limit}
\setcounter{equation}{0} \setcounter{theorem}{0} \setcounter{lemma}{0}\setcounter{proposition}{0}\setcounter{remark}{0}\setcounter{definition}{0}\setcounter{hypothesis}{0}

\section{Introduction}
Peridynamics, introduced by Silling in 2000, \cite{Silling1} is a  nonlocal formulation of continuum mechanics expressed in terms of regular elastic potentials. The theory is formulated in terms of displacement differences as opposed to spatial derivatives of the displacement field.  These features provide the flexibility to simultaneously simulate kinematics involving both smooth deformations and defect evolution.  Numerical simulations based on peridynamic modeling exhibit the formation and evolution of sharp interfaces associated with defects
and  fracture  \cite{SillingBobaru}, \cite{SillingAscari2},  \cite{Bobaru1}, \cite{WecknerAbeyaratne}, and \cite{BhattacharyaDyal}. These aspects are exploited in the peridynamic scheme for dynamic fracture simulation where the crack path is determined as part of the solution \cite{SillingAscari3}, \cite{Bobaru2}. This type of solution is distinct from the classical setting where the crack path is specified a priori see, \cite{Freund}.

We consider peridynamic formulations with unstable constitutive laws that soften beyond a critical stretch.
Here we discover new quantitative and qualitative information that is extracted from the peridynamic formulation using scaling arguments and by passing to a distinguished small horizon limit. In this limit the dynamics correspond to the simultaneous evolution of elastic displacement and fracture.  The displacement fields are shown to satisfy the wave equation. The wave equation provides the dynamic coupling between elastic waves and the evolving fracture path inside the media. The limit evolutions have bounded energy expressed in terms of the bulk and surface energies of brittle fracture mechanics. They also satisfy an energy inequality expressed in terms of the kinetic energy of the motion together with the bulk elastic energy and a Griffith surface energy. The elastic moduli and energy release rate have explicit formulas given in terms of the  moments of the peridynamic influence function and the peridynamic potential energy. These explicit formulas provide a rigorous means to calibrate the nonlinear potentials of peridynamics with  experimentally measured values of elastic constant, wave speed and critical energy release rate.

To present the ideas we focus on antiplane shear problems posed over a bounded convex domain $D\subset\mathbb{R}^2$.  The antiplane displacement transverse to $D$ is written $u(t,x)$. In the peridynamic formulation one considers pairs of points $x$, $x'$ in $\mathbb{R}^2$ and the relative displacement $\eta(t,x)=u(t,x')-u(t,x)$.  The family of points $x'$ that interact with $x$ is confined to a neighborhood $\mathcal{H}_\epsilon(x)$ of $x$ of diameter $2\epsilon$. Here $\epsilon$ is the horizon for the peridynamic interaction and $\mathcal{H}_\epsilon(x)$ is a disk of radius $\epsilon$ centered at $x$. The peridynamic influence function is defined inside $\mathcal{H}_\epsilon(x)$ and is written $J(\frac{|x'-x|}{\epsilon})$, with $M>J(r)\geq 0$ for $0\leq r\leq 1$ and zero outside.  For points $x$ residing outside $D$ but within a fixed distance $\alpha>\epsilon$ from $D$ we set the displacement $u(t,x)=0$ for all $0\leq t\leq T$; this gives nonlocal  boundary conditions of Dirichlet type for this problem \cite{DuGunzbergerlehoucqmengesha}. The domain containing the set of points $x$ for which $dist(x,D)<\alpha$ is denoted by $D_\alpha$.

\begin{figure}[htbp]
   \centering
   \includegraphics[width=0.5\textwidth]{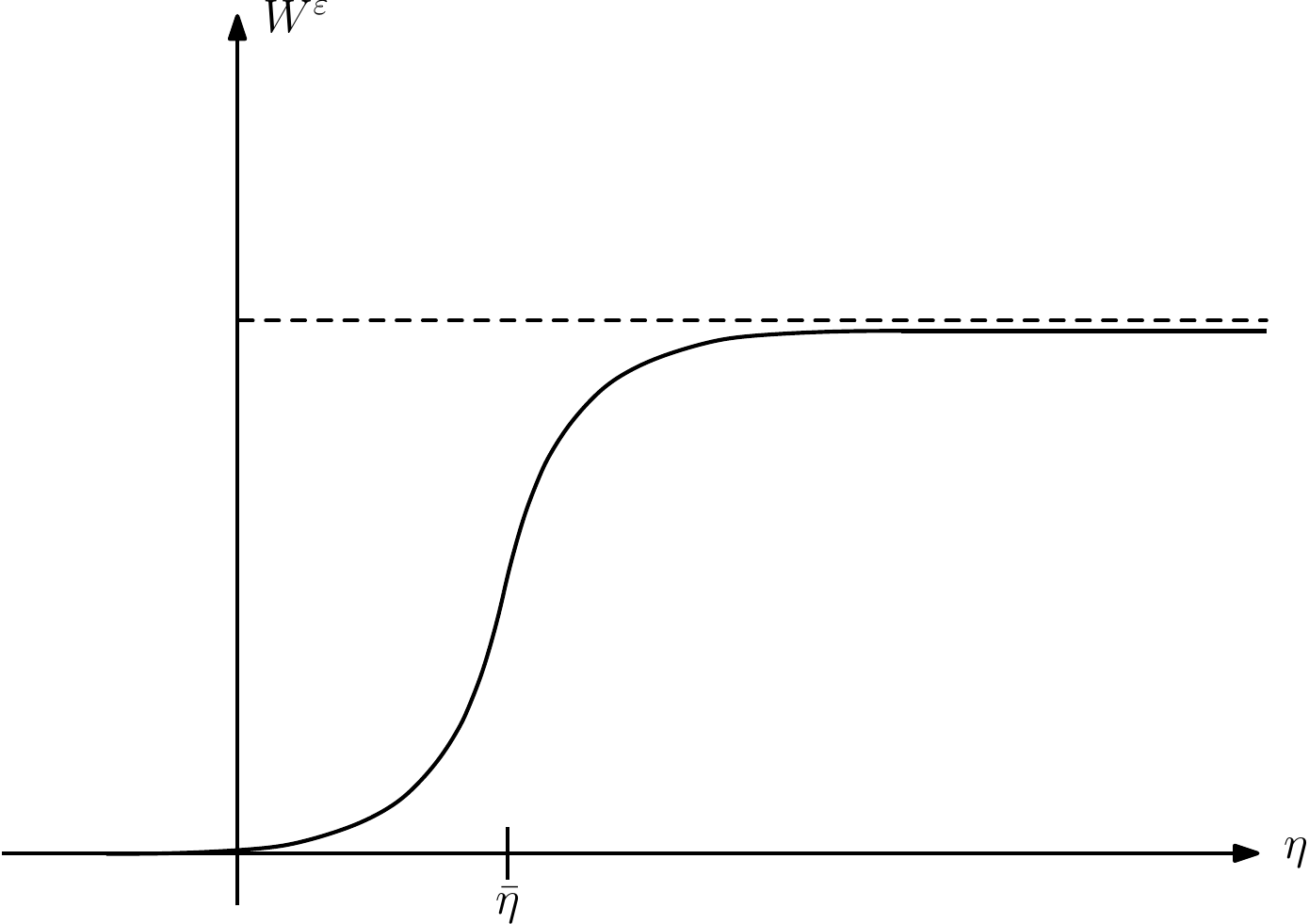}
   \caption{{\bf Convex-concave potential.}}
   \label{ConvexConcave}
\end{figure}

\begin{figure}[htbp]
   \centering
   \includegraphics[width=0.5\textwidth]{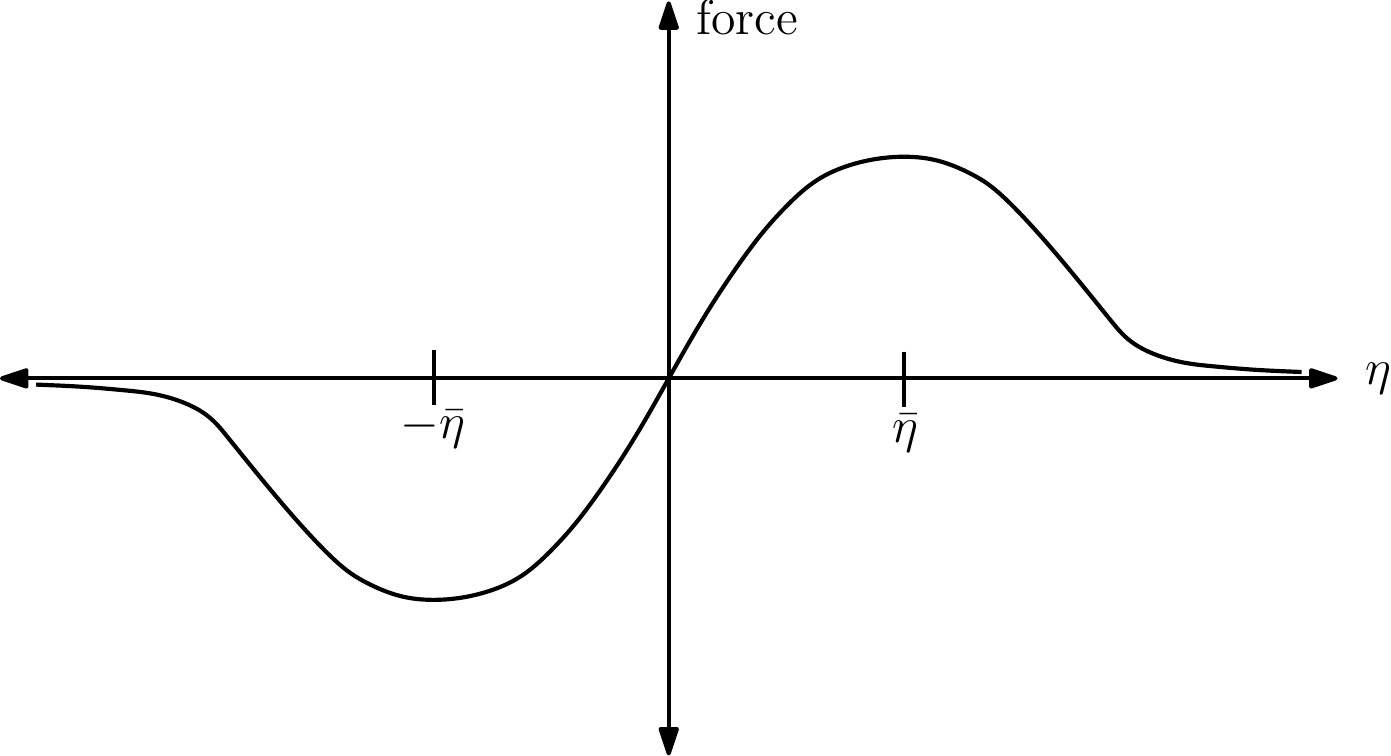}
   \caption{{\bf Unstable force versus displacement law.}}
   \label{SofteningBond}
\end{figure}

In this paper we are interested in the small horizon  limit $\epsilon\rightarrow 0$ and we make the change of variable $x'=x+\epsilon\xi$, where $\xi$ belongs to the unit disk $\mathcal{H}_1(0)$ centered at the origin. 
The peridynamic potential energy density is a function of $x'-x=\epsilon\xi$ and $\eta(x)$ and we consider the family of regular peridynamic potentials parameterized by $\epsilon$ given by
\begin{eqnarray}
W^\epsilon(\eta(x),\epsilon\xi)=\frac{1}{\epsilon^3}J\left(|\xi|\right)f\left(\frac{|\eta(x)|^2}{\epsilon|\xi|}\right).
\label{potentialdensity1}
\end{eqnarray}
The potential functions $f:[0,\infty)\rightarrow\mathbb{R}$ considered here are positive, smooth and concave with the properties
\begin{eqnarray}
\lim_{r\rightarrow 0^+}\frac{f(r)}{r}=f'(0)>0,&&\lim_{r\rightarrow\infty}f(r)=f_\infty <\infty.
\label{properties}
\end{eqnarray}
The potentials $W^\epsilon$ can be thought of as smoothed out versions of potentials used to describe the peridynamic bond stretch models introduced in \cite{Silling1}, \cite{SillingAscari2}. 
This class of potential energies is convex - concave in the relative displacement $\eta\rightarrow W^\epsilon(\eta,\epsilon\xi)$ with infection point $\overline{\eta}$ see, Figure \ref{ConvexConcave}. This delivers the constitutive relation 
\begin{eqnarray}
{\rm force}=\partial_\eta W^\epsilon(\eta(x),\epsilon\xi)
\label{forcestate}
\end{eqnarray}
see, Figure \ref{SofteningBond}. Here the magnitude of the force increases with the magnitude of the relative displacement up to the critical value $\overline{\eta}$ after which it decreases.

The peridynamic potential energy is obtained by integrating the energy density over the neighborhood $\mathcal{H}_\epsilon(x)$ and is given 
in terms of the rescaled coordinates by
\begin{eqnarray}
\epsilon^2\int_{\mathcal{H}_1(0)}W^\epsilon(\eta(x),\epsilon\xi)\,d\xi.
\label{potentialenergy}
\end{eqnarray}
The total strain energy of the displacement is given by
\begin{eqnarray}
PD^\epsilon(u)=\epsilon^2\int_{D}\int_{\mathcal{H}_1(0)}W^\epsilon(\eta(x),\epsilon\xi)\,d\xi\,dx.
\label{TotalStrain}
\end{eqnarray}

The peridynamic evolution is described by the Lagrangian
\begin{eqnarray}
L^\epsilon(u(t),\partial_t u(t),t)=K(\partial_t u(t))-PD^\epsilon(u(t))+U(u(t)),
\label{Lagrangian}
\end{eqnarray}
with 
\begin{eqnarray}
K(\partial_t u(t))&=&\frac{1}{2}\int_{D}\rho|\partial_t u(t,x)|^2\,dx, \hbox{ and }\nonumber\\
U(u(t))&=&\int_{D}b(t,x) u(t,x)\,dx,
\label{Components}
\end{eqnarray}
where $\rho$ is the mass density of the material and $b(t,x)$ is the body force density. 
The  initial conditions $u^\epsilon(0,x)=u_0(x)$ and $u^\epsilon_t(0,x)=v_0(x)$ are prescribed and the action integral for the peridynamic evolution is
\begin{eqnarray}
I^\epsilon(u)=\int_0^TL^\epsilon(u(t),\partial_t u(t),t)\,dt.
\label{Action}
\end{eqnarray}
For a given a unit vector $n$ and $h$ in $\mathbb{R}$ we define the difference quotient of $\psi(x)$ 
by
\begin{equation}
D_n^{h}\psi(x)=\left\{\begin{array}{cl}\frac{\psi(x+hn)-\psi(x)}{h}& \hbox{, if $h\not=0$}\\
0&\hbox{, if $h=0$.}
\end{array} \right.
\label{diffquotient}
\end{equation}
Writing $e=\frac{\xi}{|\xi|}$ and $\xi=e|\xi|$, we set 
\begin{eqnarray}
D_e^{\epsilon|\xi|}\psi(x)&=&\frac{\psi(x+\epsilon\xi)-\psi(x)}{\epsilon|\xi|} \label{diffp}\\
D_{-e}^{\epsilon|\xi|}\psi(x)&=&\frac{\psi(x-\epsilon\xi)-\psi(x)}{\epsilon|\xi|}.\label{diffm}
\end{eqnarray}
The Euler Lagrange Equation for this system delivers the peridynamic equation of motion given by
\begin{eqnarray}
\rho u^\epsilon_{tt}&=&-\nabla PD^\epsilon(u^\epsilon)+b
\label{stationary}
\end{eqnarray}
where
\begin{eqnarray}
\nabla PD^\epsilon(u^\epsilon)=\epsilon^3\int_{\mathcal{H}_1(0)}D_{-e}^{\epsilon|\xi|}\partial_\eta W^\epsilon(\eta(x),\epsilon\xi)|\xi|\,d\xi.
\label{GradPD}
\end{eqnarray}

In what follows we examine
the family of peridynamic deformations $\{u^\epsilon(t,x)\}_{\epsilon>0}$ defined for suitable initial data $u_0,v_0$ and investigate
the the behavior of the limiting deformation field $u^0(t,x)$. To do this we describe peridynamic deformations as trajectories in function space. 
The nonlocal Dirichlet boundary conditions are incorporated into the function space by defining the class of functions $L^2_0(D)$ that are square integrable over $D$ and zero on $D_\alpha\setminus D$.  \footnote{For $A\supset B$,  we denote $A\setminus B=A\cap B^c$, where $B^c$ is the complement of $B$ in $A$.} In this context we view peridynamic evolutions as functions of time taking values in the space  $L^2_0(D)$.  It follows from the evolution equation \eqref{stationary} that $u^\epsilon(t,x)$ is twice differentiable in time taking values in $L^2_0(D)$. This space of functions is denoted by $C^2([0,T];L^2_0(D))$ see, e.g., \cite{Evans}. The initial value problem for the peridynamic evolution  \eqref{stationary}
is shown to be well posed on $C^2([0,T];L^2_0(D))$ see, section \ref{PE}.
We apply a scaling analysis to show that the peridynamic evolutions $u^\epsilon(t,x)$  approach a limit evolution $u^0$ in the $\epsilon\rightarrow$ limit. The limit evolution $u^0(t,x)$ is shown to have bounded elastic energy in the sense of fracture mechanics for a wide class of initial conditions.  The limit evolution satisfies an energy inequality expressed in terms of the kinetic energy of the motion together with a linear elastic energy in terms of the antiplane strain $\nabla u^0$ and a Griffith surface energy associated with the evolving jump set $S_{u^0(t)}$ of $u^0(t,x)$ see, section \ref{DBF}.

Motivated by the approach given in \cite{SillingWecknerBobaru} we investigate the effect of bond instability on the nucleation of fracture inside a peridynamic body.  We consider a generic peridynamic neighborhood $\mathcal{H}_\epsilon(x)$ of radius $\epsilon$ about the point x. For points $x'$ inside $\mathcal{H}_\epsilon(x)$ we say that the material between $x$ and $x'$  (a bond) is critically stretched if $|\eta^\epsilon(x)|=|u^\epsilon(x')-u^{\epsilon}(x)|>\overline{\eta}$, otherwise the relative displacement is called subcritical. A linear stability analysis is given that identifies necessary conditions for fracture nucleation inside $\mathcal{H}_\epsilon(x)$. These conditions are directly linked to the appearance of subsets with nonzero area fraction containing critically stretched bonds.  The fracture nucleation condition given by Proposition \ref{Nuccriteria} implies that if the neighborhood contains a nonzero area fraction of critically stretched bonds then the neighborhood can be linearly unstable and a displacement jump can be nucleated. 
These results are presented in section \ref{stablilty}.

Motivated by this result we focus on  peridynamic neighborhoods $\mathcal{H}_{\epsilon}(x)$  that contain critically stretched bonds over an area fraction  larger than $\sqrt{\epsilon}$. These neighborhoods are referred to as unstable neighborhoods. Under this definition unstable neighborhoods have the potential to nucleate jump discontinuities. We apply this definition to identify a concentration of bond instability as the peridynamic horizon approaches zero. To present the idea we consider the collection of centroids  $x$  of all the unstable neighborhoods $\mathcal{H}_\epsilon(x)$ with radii  less than a tolerance $\delta$ . This collection is denoted by the set $C_\delta$. It is shown that the area of $C_\delta$ vanishes as $\delta\rightarrow 0$ and that the collection of centroids for unstable neighborhoods concentrate onto a set $C_0$ of zero area (a set of Lebesgue measure zero) as $\delta\rightarrow 0$.   The concentration of instability with respect to peridynamic horizon is directly linked to the energy budget associated with the peridynamic motion. It is shown that for a family of peridynamic flows $\{u^\epsilon(t,x)\}_{\epsilon>0}$ all driven by the same initial conditions and body forces that the peridydnamic potential energy of each flow is bounded uniformly  in time $0\leq t\leq T$  independently of the radius of the horizon, see section \ref{GKP}. This bound is shown to force the localization see Theorem  \ref{bondunstable}.   These observations are presented in section \ref{localglobal} and established in section \ref{proofbondunstable}. 

We apply these observations and adopt the hypothesis that the jump set of the limit evolution $S_{u^0}$  and the concentration set $C_0$ are one and the same.   We employ a scaling analysis  to the peridynamic equation of motion to discover that the limit evolution $u^0(t,x)$ satisfies the wave equation at every point where $\nabla u^0$ is defined see Theorem \ref{waveequation}.  The wave equation provides the dynamic coupling between elastic waves $u^0(t,x)$  and the evolving fracture path $S_{u^0(t)}$ inside the media. It is important to point out that the limiting dynamic fracture evolution described here follows from scaling arguments and on passing to a distinguished limit in the peridynamic formulation. 
These results are presented in sections \ref{DBF} and \ref{STAY}.  The mathematical tool set appropriate for extracting the limit behavior from this class of peridynamic models is based on $\Gamma$-convergence and comes from the literature associated with the analysis of the Mumford Shah functional and free discontinuity problems see, \cite{Gobbino1}, \cite{Gobbino2}, and \cite{Gobbino3}.

Other related recent work focuses on passing to the small horizon limit for linear peridynamic formulations; establishing a link between linear elasticity and peridynamics see \cite{EmmrichWeckner}, \cite{SillingLehoucq}.  A global stability criterion for describing the nucleation of phase transformations within the peridynamic setting is developed in \cite{BhattacharyaDyal}.

We note  that there is a vast  literature on fracture modeling and a complete survey is beyond the scope of this paper. Instead we point out  recent proposals for computing crack propagation in dynamic and quasi static settings.  Approaches using a phase field for the damage set and a linear elastic field, to represent crack propagation have been proposed and developed in \cite{BourdinLarsenRichardson}, \cite{LarsenOrtiner}, and \cite{Hughes}.  Wave equations for fields inside domains with evolving cracks are posed in \cite{LarsenDalMaso} and variational aspects of sharp interface models are discussed in \cite{Larsen}.  For quasi static problems variational phase field methods are developed in the pioneering work of \cite{BourdinFrancfortMarigo}, \cite{FrancfortMarigo}. More recently a two field method using eigen-deformations for the fracture field is developed for quasi static problems in \cite{Ortiz}. Alternative nonlocal formulations have been developed for quasi static crack propagation in \cite{Braides}, \cite{Buttazzo}.

\setcounter{equation}{0} \setcounter{theorem}{0} \setcounter{lemma}{0}\setcounter{proposition}{0}\setcounter{remark}{0}\setcounter{definition}{0}\setcounter{hypothesis}{0}

\section{Peridynamic evolution}
We begin this section by introducing a suitable class of initial conditions appropriate for describing the  evolution
of deformations that can have smooth variation as well as jumps. Here we will choose initial conditions with bounded elastic energy in the sense of fracture mechanics. We show that well posed peridynamic evolutions  exist for this class of initial data. These peridynamic evolutions  satisfy an energy balance between potential and kinetic energy at each time during the deformation. Next we develop  a necessary criterion for fracture initiation inside a peridynamic neighborhood. Here fracture initiation is defined to be the nucleation of a jump in the displacement inside a peridynamic neighborhood. We develop a criterion for the orientation of the nucleated crack based upon the notion of the  {\em most unstable} direction. The approach taken here is consistent with the analysis of crack nucleation developed in \cite{SillingWecknerBobaru}.  We conclude with a discussion  of the localization of instability in the limit of vanishing  peridynamic horizon. 

\subsection{Initial conditions and motivation}
\label{IBF}
Our choice of initial conditions is motivated by  Linear Elastic Fracture Mechanics (LEFM). Here we investigate Mode III fracture in the context of antiplane shear. The initial condition is specified by a crack set $K$ and displacement $u_0$. The gradient  $\nabla u_0$ is defined off the crack set and the displacement $u_0$ can suffer jumps across $K$.  Griffith's theory of brittle fracture
asserts that the energy necessary to produce a crack $K$ is proportional to the crack length $\ell$. For LEFM the total energy associated with bulk elastic and surface  energy is given by
\begin{eqnarray}
\int_{D}\mu|\nabla u_0|^2\,dx+\mathcal{G}_c\ell,
\label{Gcrackenergy}
\end{eqnarray}
where $\mu$ is the shear modulus and $\mathcal{G}_c$ is the critical energy release rate for the material.
In what follows we chose initial conditions associated with bounded LEFM elastic energy.

In order to pass to the small horizon limit of peridynamics and to understand the elastic energy associated with this limit we cast the problem in a functional analytic context. The function space used in the mathematical formulation of free discontinuity problems including fracture is the space ${SBV}$ developed in \cite{DeGiorgiAmbrosio} see also, \cite{AmbrosioBrades}, \cite{Ambrosio}.   Functions in this space belong to  $L^1(D)$  and are approximately continuous almost everywhere.
Here we recall that  points $x$ of approximate continuity for the function $u_0$ satisfy
\begin{eqnarray}
\lim_{r\searrow 0}\frac{1}{\pi r^2}\int_{B(x,r)}\,|u_0(y)-u_0(x)|\,dy=0,
\label{approx}
\end{eqnarray}
where $B(x,r)$ is the ball of radius $r$ centered at $x$.
The  discontinuity set $S_{u_0}$  for elements of  $SBV$ are characterized by at most a countable number of smooth rectifiable curves across which $u_0$ has a jump discontinuity. Here the notion of arc length corresponds to the one dimensional Hausdorff  measure of $S_{u_0}$ and is denoted by $\mathcal{H}^1(S_{u_0})$.   We choose an orientation and define the unit normal $\nu$ to the jump set $S_{u_0}$. 
For points $x$ belonging to  the jump set we denote the intersection of $B(x,r)$ with the half spaces $(y-x)\cdot \nu<0$ and $(y-x)\cdot \nu>0$ by $B^-(x,r)$ and $B^+(x,r)$ respectively. The left and right limits of the function  $u_0$ for a point on the jump set are denoted by $u^-_0$, $u^+_0$ and satisfy the identities
\begin{eqnarray}
\lim_{r\searrow 0}\frac{1}{\pi r^2}\int_{B^-(x,r)}\,|u_0(y)-u^-_0(x)|\,dy=0, &\hbox{     }&\lim_{r\searrow 0}\frac{1}{\pi r^2}\int_{B^+(x,r)}\,|u_0(y)-u^+_0(x)|\,dy=0.
\label{approxjump}
\end{eqnarray}
The approximate gradient denoted by $\nabla u_0$ of an SBV function is defined almost everywhere on $D\setminus S_{u_0}$ and satisfies
\begin{eqnarray}
\lim_{r\searrow 0}\frac{1}{\pi r^2}\int_{B(x,r)}\,\frac{|u_0(y)-u_0(x)-(y-x)\cdot\nabla u_0(x)|}{r}\,dy=0.
\label{appgrad}
\end{eqnarray}
Distributional derivatives $Du_0$ of $SBV$ functions  are constructed from the approximate gradient and jump sets and satisfy
\begin{eqnarray}
\langle Du_0,\Phi \rangle=\int_{D}\nabla u_0\cdot\Phi\,dx+\int_{S_{u_0}}(u_0^+-u_0^-)\nu\cdot\Phi d\mathcal{H}^1,
\label{derivative}
\end{eqnarray}
for every continuous test function $\Phi:D \rightarrow\mathbb{R}^2$ with support on $D$. Here $d\mathcal{H}^1$ corresponds to an element of  arc length for sufficiently regular curves. Functions in SBV have distributional derivatives with bounded total variation. Detailed descriptions of the properties of SBV functions are provided in \cite{Ambrosiobook} and \cite{Buttazzo}.

Deformations of class $SBV(D)$ are easily interpreted as deformations with cracks in $D$: the crack set $K$ is identified with the jump set $S_{u_0}$  and $\nabla u_0$ represents the usual strain in the elastic part of the body outside the crack see, \cite{AmbrosioBrades}, \cite{Ambrosio}. With this in mind we take the initial displacement $u_0\in L^2_0(D)$  and require that it belong to the space  $SBV(D)$. For this choice of initial data the bulk and surface energy of LEFM is given by
\begin{eqnarray}
LEFM(u_0,D)=\int_{D}\mu|\nabla u_0|^2\,dx+\mathcal{G}_c\mathcal{H}^1(S_{u_0}).
\label{Gcrackenergyweak}
\end{eqnarray}

\begin{definition}
\label{LEFMID}
We refer to initial data  $u_0\in L^2_0(D)$, $v_0\in L^2_0(D)$ with the restriction of $u_0$ on $D$ belonging to $SBV(D)$ that satisfy
\begin{eqnarray}
LEFM(u_0,D)<\infty,  \hbox{        }\sup_{x \in D} |u_0|<\infty, \hbox{        }\sup_{x \in D} |v_0|<\infty,
\label{lefmid}
\end{eqnarray}
as LEFM initial data. 
\end{definition}

We coordinate  our choice of shear modulus $\mu$ and critical energy release rate $\mathcal{G}_{c}$  with the peridynamic potential $f$ and influence function $J$ through the relations:
\begin{eqnarray}
\mu=\pi \, f'(0)\int_{0}^1r^2J(r)dr &\hbox{          }& \mathcal{G}_c=2\pi\, f_\infty \int_{0}^1r^2J(r)dr,
\label{calibrate}
\end{eqnarray}
where $f_\infty$ is defined by \eqref{properties}. The correspondence between the shear modulus and critical energy release rate and the peridynamic quantities $f$ and $J$  follows directly from the limit analysis, see  Theorem \ref{limitflow}, section \ref{CC}, and equation  \eqref{shear}.

\subsection{Peridynamic evolutions and energy balance}
\label{PE}
We choose  the  initial data $(u_0,v_0)$ to be   LEFM initial data  and the initial crack set at $t=0$ is prescribed by $K=S_{u_0}$. There is a unique peridynamic evolution for this choice of initial data. This is stated in the following theorem.
\begin{theorem}{\rm {\bf Existence of unique solution for nonlinear peridynamics}}\\
\label{existenceuniqueness}
\noindent For LEFM initial data  $(u_0,v_0)$  and body force $b(t,x)$ in $C^1([0,T];L^2(D))$ there exists a unique peridynamic evolution $u^\epsilon(t,x)$ in $C^2([0,T];L^2_0 (D))$
taking the initial values $u^\epsilon(0,x)=u_0(x)$, $u_t^\epsilon(0,x)=v_0(x)$, and satisfying
\begin{eqnarray}
&&\rho u^\epsilon_{tt}(t,x)=-\nabla PD^\epsilon(u^\epsilon(t,x))+b(t,x), \hbox{ \rm  for $0<t\leq T$ and $x$ in $D$.}
\label{stationaryagain}
\end{eqnarray}
\end{theorem}
\noindent This theorem follows from the Lipschitz continuity of $\nabla PD^\epsilon$ and is established in section \ref{EE}.

Multiplying both sides of \eqref{stationaryagain} by $u^\epsilon_t$ delivers the identity
\begin{eqnarray}
\partial_t\left\{\frac{\rho}{2}\Vert u^\epsilon_t\Vert^2_{L^2(D)}+PD^\epsilon(u^\epsilon)\right\}=\int_{D}bu^\epsilon_t\,dx
\label{differential}
\end{eqnarray}
and integration over time from $0$ to $t$ delivers the energy balance associated with the peridynamic evolution given by
\begin{theorem}{\rm \bf Energy balance}\\
\label{Ebalance}
\begin{eqnarray}
\mathcal{EPD}^\epsilon(t,u^\epsilon(t))=\mathcal{EPD}^\epsilon(0,u^\epsilon(0))-\int_0^t\int_{D} b_t(\tau)u^\epsilon(\tau)\,dx\,d\tau,\label{BalanceEnergy}
\end{eqnarray}
where
\begin{eqnarray}
\mathcal{EPD}^\epsilon(t,u^\epsilon(t))=\frac{\rho}{2}\Vert u^\epsilon_t(t)\Vert^2_{L^2(D)}+PD^\epsilon(u^\epsilon(t))-\int_{D}b(t)u^\epsilon(t)\,dx
\label{energyt}
\end{eqnarray}
and
\begin{eqnarray}
\mathcal{EPD}^\epsilon(0,u^\epsilon(0))=\frac{\rho}{2}\Vert v_0\Vert^2_{L^2(D)}+PD^\epsilon(u_0)-\int_{D}b(0)u_0\,dx.
\label{energy0}
\end{eqnarray}
\end{theorem}

\subsection{Instability and fracture initiation}
\label{stablilty}
In this section we present a fracture initiation condition that arises from the unstable peridynamic constitutive law relating relative displacement  to force.
This type of fracture nucleation condition has been  identified for peridynamic evolutions in \cite{SillingWecknerBobaru}. Here we investigate the nucleation criteria for the case at hand and provide an additional condition for the  most unstable direction along which the 
crack can nucleate. We introduce a jump perturbation at $x$ associated with a direction $\nu$ on the unit circle. Set $E_{\nu}^+(x)=\{y:\,(y-x)\cdot\nu^\perp>0\}$ and  $E_{\nu}^-(x)=\{y:\,(y-x)\cdot\nu^\perp\leq0\}$ and introduce the local coordinate basis at $x$ given by $\nu$ and $\nu^\perp$. Consider a time independent body force density $b$ and a smooth equilibrium solution $u$ of \eqref{stationary}. We now perturb $u$ by adding a function $\delta_\nu$ with a jump discontinuity of height $\delta$ across the line $\{y\in\mathcal{H}_\epsilon(x);\,(y-x)\cdot\nu^\perp=0\}$ that is piecewise constant in $\mathcal{H}_\epsilon(x)$ and $\delta_\nu=\delta$ for points in $E_\nu^+\cap\mathcal{H}_1(0)$ and $\delta_\nu=0$ for points in $E_\nu^-\cap\mathcal{H}_1(0)$. Here the direction $\nu$ points along the direction of the discontinuity and $\nu^\perp$ is the normal to the line of discontinuity. We write  $u^p=u+\delta_\nu$ and apply the ansatz
\begin{eqnarray}
\rho u^p_{tt}&=&-\nabla PD^\epsilon(u^p)+b.\label{stationarydiff1}
\end{eqnarray}
We regard $\delta$ as a small perturbation and expand the integrand of $\nabla PD^\epsilon(u^p)$ in a Taylor series to recover the linearized evolution equation for the jump $\delta$ at $x$ across the line with normal $\nu^\perp$. The evolution equation is given by
\begin{eqnarray}
\rho\delta_{tt}=\mathcal{A}_\nu\delta,
\label{pertevolution}
\end{eqnarray}
where
\begin{eqnarray}
\mathcal{A}_\nu&=&-\frac{1}{2}\left\{\int_{\mathcal{H}_1(0)\cap E_\nu^+(0)}\epsilon^2\partial^2_{\eta}W^\epsilon(u(x+\epsilon\xi)-u(x),\epsilon\xi)d\xi\right.\nonumber\\
&&\left.+\int_{\mathcal{H}_1(0)\cap E_\nu^-(0)}\epsilon^2\partial^2_{\eta}W^\epsilon(u(x)-u(x-\epsilon\xi),\epsilon\xi)d\xi\right\},
\label{instabilitymatrix}
\end{eqnarray}
here $E_{\nu}^+(0)=\{\xi:\,\xi\cdot\nu^\perp>0\}$ and  $E_{\nu}^-(0)=\{\xi:\,\xi\cdot\nu^\perp\leq0\}$.
Calculation shows that
\begin{eqnarray}
\partial^2_{\eta}W^\epsilon(\eta,\epsilon\xi)=\frac{2}{\epsilon^4|\xi|}J(|\xi|)\left(f'\left(\frac{\eta^2}{\epsilon|\xi|}\right)+2f''\left(\frac{\eta^2}{\epsilon|\xi|}\right)\frac{\eta^2}{\epsilon|\xi|}\right),\label{expand}
\end{eqnarray}
where $f'(\eta^2/\epsilon|\xi|)>0$, $f''(\eta^2/\epsilon|\xi|)<0$ and the critical value $\overline{\eta}$ is the root of $\partial_\eta^2W^\epsilon(\eta,\epsilon\xi)=0$ with $\partial_\eta^2W^\epsilon(\eta,\epsilon\xi)>0$ for $|\eta|<\overline{\eta}$ and $\partial_\eta^2W^\epsilon(\eta,\epsilon\xi)<0$ for $|\eta|>\overline{\eta}$. Here $\overline{\eta}=\sqrt{\epsilon|\xi|}\overline{r}$ where $\overline{r}$ is the inflection point for the function $r:\rightarrow f(r^2)$. For $\mathcal{A}_\nu>0$ the jump can grow exponentially. It is evident that this can occur if there are critically stretched bonds, $|\eta|>\overline{\eta}$, inside the neighborhood. 
We summarize these results in the following.
\begin{proposition}{\em \bf Facture nucleation condition}\\
\label{Nuccriteria}
Given a point $x$ and a direction $\nu$ a condition for crack nucleation at $x$ along direction $\nu$ is $\mathcal{A}_\nu>0$.
The  directions $\nu^*$ along which cracks most likely grow are the most unstable ones which satisfy the condition
\begin{eqnarray}
\mathcal{A}_{\nu^*}=\max_{\nu}\mathcal{A}_\nu>0.
\label{bestdirction}
\end{eqnarray}
\end{proposition}

\subsection{Concentration of fracture nucleation sites  in the small horizon limit }
\label{localglobal}

Here we present results that show that peridynamic neighborhoods likley to  nucleate  jump sets  become concentrated in the small horizon limit.   The discussion focuses on the basic unit of peridynamic interaction: the peridynamic neighborhoods $\mathcal{H}_\epsilon(x)$ of diameter $\epsilon>0$ with centroids $x\in D$. In what follows we denote the two dimensional Lebesgue measure (area)  of a set $S$ by $|S|$. Here we investigate the family of peridynamic evolutions $u^\epsilon(t,x)$  at a fixed time $t$.

Consider a prototypical  neighborhood $\mathcal{H}_\epsilon(x)$. 
The collection of points $y$ inside $\mathcal{H}_\epsilon(x)$ for which the relative displacement is beyond critical, i.e., $|u^\epsilon(t,y)-u^\epsilon(t,x)|>\overline{\eta}$ is called the unstable subset of $\mathcal{H}_\epsilon(x)$ and is written as
\begin{eqnarray}
\left\{y\hbox{ in }\mathcal{H}_{\epsilon}(x):\, |u^\epsilon(t,y)-u^\epsilon(t,x)|>\overline{\eta}\right\},
\label{unstable}
\end{eqnarray}
where 
$\overline{\eta}=\sqrt{|y-x|}\overline{r}$,  and $\overline{r}$ is the inflection point for the map $r:\rightarrow f(r^2)$. 
The weighted area fraction of the neighborhood $\mathcal{H}_\epsilon(x)$ occupied by the unstable subset is denoted by
\begin{eqnarray}
P(\left\{y\hbox{ in }\mathcal{H}_{\epsilon}(x):\, |u^\epsilon(t,y)-u^\epsilon(t,x)|>\overline{\eta}\right\}).
\label{weightarea}
\end{eqnarray}
Here $P$ is defined in terms of  the indicator function $\chi^{+,\epsilon}(x,y)$ for the unstable subset with,
$\chi^{+,\epsilon}(x,y)=1$ for $y$ in the unstable subset and $0$ otherwise, and 
\begin{eqnarray}
P(\left\{y\hbox{ in }\mathcal{H}_{\epsilon}(x):\, |u^\epsilon(t,y)-u^\epsilon(t,x)|>\overline{\eta}\right\})=\frac{1}{\epsilon^2 m}\int_{\mathcal{H}_\epsilon(x)}\,\chi^{+,\epsilon}(x,y)\left\vert\frac{y-x}{\epsilon}\right\vert J(\left\vert\frac{y-x}{\epsilon}\right\vert)\,dy,
\label{weight}
\end{eqnarray}
where the normalization constant $m=\int_{\mathcal{H}_1(0)}|\xi|J(|\xi|)\,d\xi$ is chosen so that $P(\mathcal{H}_\epsilon(x))=1$.

\begin{definition}
\label{UnstabeL}
The neighborhood $\mathcal{H}_{\epsilon}(x)$ is said to be unstable at time $t$ if
\begin{eqnarray}
\sqrt{\epsilon}<P\left(\left\{y\hbox{ in }\mathcal{H}_{\epsilon}(x):\, |u^\epsilon(t,y)-u^\epsilon(t,x)|>\overline{\eta}\right\}\right).
\label{unstablee}
\end{eqnarray}
\end{definition}

To proceed choose a fixed length scale $\delta>0$  and consider a family of radii  $\epsilon_j=\frac{1}{2^j}$, $j=1,\ldots$ and  the collection of neighborhoods $\mathcal{H}_{\epsilon_j}(x)$ with centroids $x$ in the reference domain $D$. The set of  centroids associated with unstable neighborhoods for  $\epsilon_j<\delta$ at time $t$ is denoted by $C_{\delta,t}$. This set is expressed as
\begin{eqnarray}
C_{\delta,t}=\left\{x\in D;\, \exists \, \epsilon_j<\delta\hbox{ for which } P\left(\left\{y\hbox{ in }\mathcal{H}_{\epsilon_j}(x):\, |u^{\epsilon_j}(t,y)-u^{\epsilon_j}(t,x)|>\overline{\eta}\right\}\right)>\sqrt{\epsilon_j}\right\}.
\label{unstabelll}
\end{eqnarray}
Here $C_{\delta,t}\subset C_{\delta',t}$ for $\delta<\delta'$.
Let $C_{0,t}={\cap}_{\scriptscriptstyle{0<\delta}} C_{\delta,t}$ denote the concentration set for the set of centroids associated with unstable neighborhoods.
We now state a theorem on the localization of bond instability as the peridynamic horizon shrinks to zero.
\begin{theorem}Localization of bond instability in the small horizon limit.
\label{bondunstable}\\
The collection of centroids $C_{\delta,t}$ for unstable neighborhoods  is decreasing  as $\delta\rightarrow 0$ and there is a positive constant $C$ indpendent of $t$ and $\delta$ for which
\begin{eqnarray}
|C_{\delta,t}|\leq C\sqrt{\delta}, \hbox{  for,  } 0\leq t\leq T.
\label{limdelta}
\end{eqnarray}
Moreover $C_{\delta,t}$ concentrate on the set $C_{0,t}$, where $C_{0,t}$ is a set of Lebesgue measure zero, i.e.,
\begin{eqnarray}
\lim_{\delta\rightarrow 0}|C_{\delta,t}|=|C_{0,t}|=0.
\label{limzero}
\end{eqnarray}
\end{theorem}
\noindent Theorem \ref{bondunstable} is established in section \ref{proofbondunstable}. The localization of instability with respect to horizon is directly linked to the energy budget associated with the peridynamic motion. It is shown that for a family of peridynamic flows $\{u^\epsilon(t,x)\}_{\epsilon>0}$ all driven by the same initial conditions and body forces that the peridydnamic potential energy of each flow is bounded uniformly  in time $0\leq t\leq T$  independently of the radius of the horizon, see section \ref{GKP}. This bound forces the localization as shown in section \ref{proofbondunstable}. 

\setcounter{equation}{0} \setcounter{theorem}{0} \setcounter{lemma}{0}\setcounter{proposition}{0}\setcounter{remark}{0}
\setcounter{definition}{0}\setcounter{hypothesis}{0}

\section{The small horizon, sharp interface limit}
In this section we identify the $\epsilon\searrow 0$ limit of the solutions $u^\epsilon$ to the peridynamic initial value problem with LEFM initial data. A limit evolution $u^0(t,x)$ is identified that:
\begin{itemize}
\item Has a uniformly bounded bulk elastic energy and a Griffith surface energy  associated with fracture mechanics for $0\leq t\leq T$.
\item Satisfies an energy inequality involving the kinetic energy of the motion together with the bulk elastic and sufrace energy associated with  fracture mechanics for  $0\leq t\leq T$.
\item Satisfies the wave equation for  $0\leq t\leq T$.
\end{itemize}

\subsection{Convergence of peridynamics to sharp interface dynamics associated with brittle fracture}
\label{DBF}

We consider the family of solutions $u^{\epsilon_k}$ to the peridynamic initial value problem with LEFM initial data for a sequence $\epsilon_{k}$, $k=1,2,\ldots$.  We shall see that we can pass to the limit $\epsilon_{k}\searrow 0$ to identify a limit evolution $u^0(t,x)$ for $0\leq t\leq T$.  The limit flow is found to  have an approximate gradient $\nabla u^0(t,x)$ almost everywhere in $D$  and the jump set $S_{u^0(t)}$  is the countable union of rectifiable arcs. 
Moreover the limit evolutions $u^0(t,x)$ have bounded energy in the sense of fracture mechanics over $0\leq t\leq T$. We begin by making the following hypothesis.

\begin{hypothesis}
\label{remarkone}
We suppose that the magnitude  of the deformations do not become infinite for $0\leq t\leq T$, i.e.,
\begin{eqnarray}
\sup_{k}\Vert u^{\epsilon_k}(t)\Vert_{L^\infty(D)}<\infty,
\label{max}
\end{eqnarray}
for $0\leq t\leq T$.  This hypothesis  is consistent with the  bounds on the kinetic  energy for peridynamic evolution given in Theorem \ref{Gronwall} of Section \ref{GKP} and is also motivated by simulations carried out in the peridynamic literature see, for example \cite{Bobaru2}, \cite{SillingAscari2}. 
\end{hypothesis}

\begin{theorem} 
\label{limitflow}
{\rm\bf Limit evolution with bounded LEFM energy.}\\
Suppose Hypothesis  \ref{remarkone} holds true then there exists a subsequence of peridynamic evolutions $u^{\epsilon_k}$ with LEFM initial data that converge as trajectories in  $C([0,T];L^2_0(D))$ to $u^0(t,x)$ in $C([0,T];L_0^2(D))$. The limit flow has an approximate gradient $\nabla u^0(t,x)$ almost everywhere in $D$  and the jump set $S_{u^0(t)}$  is the countable union of rectifiable arcs. Furthermore there exists a constant $C$ depending only on $T$  bounding the LEFM energy of the limit flow, i.e., 
\begin{eqnarray}
\mu\int_{D}|\nabla u^0(t,x)|^2\,dx+\mathcal{G}_c\mathcal{H}^1(S_{u^0(t)})\leq C
\label{LEFMbound}
\end{eqnarray}
for $0\leq t\leq T$.
The relations between the peridynamic potential $f$, influence function $J$, shear modulus $\mu$,  and critical energy release rate $\mathcal{G}_c$  are given by \eqref{calibrate}.
\end{theorem}
\noindent Theorem \ref{limitflow} is established using Gronwall's inequality see, section \ref{GKP} and  the $\Gamma$-- convergence associated with peridynamic energies see, section \ref{CC}. The proof of Theorem \ref{limitflow} is given in section \ref{CC}.

We now present an energy inequality for the limit evolution. We denote the LEFM energy for the limit evolution $u^0(t)$  at time $t$ as
\begin{eqnarray}
LEFM(u^0(t),D)=\mu\int_{D}|\nabla u^0(t)|^2\,dx+\mathcal{G}_c\mathcal{H}^1(S_{u^0(t)})
\label{LEFMatt}
\end{eqnarray}
and the LEFM energy for the initial data is written
\begin{eqnarray}
LEFM(u_0,D)=\mu\int_{D}|\nabla u_0|^2\,dx+\mathcal{G}_c\mathcal{H}^1(S_{u_0}).
\label{LEFMinit}
\end{eqnarray}
The sum of energy and work for the deformation $u^0$ at time $t$ is written
\begin{eqnarray}
\mathcal{GF}(u^0(t),D)=\frac{\rho}{2}\Vert u_t^0(t)\Vert^2_{L^2(D)}+LEFM(u^0(t),D)-\int_{D}b(t)u^0(t)\,dx.
\label{sumt}
\end{eqnarray}
The sum of energy and work for the initial data $u_0,v_0$ is written
\begin{eqnarray}
\mathcal{GF}(u_0,D)=\frac{\rho}{2}\Vert v_0\Vert^2_{L^2(D)}+LEFM(u_0,D)-\int_{D}b(0)u_0\,dx.
\label{sumt}
\end{eqnarray}
The energy inequality for the limit evolution $u^0$ is given in the following theorem.
\begin{theorem} {\rm \bf Energy Inequality.}\\
\label{energyinequality}
For $0\leq t\leq T$,
\begin{eqnarray}
\mathcal{GF}(u^0(t),D)\leq\mathcal{GF}(u_0,D)-\int_0^t\int_{D} b_t(\tau) u^0(\tau)\,dx\,dt.
\label{enegineq}
\end{eqnarray}
\end{theorem}
\noindent The proof of Theorem \ref{energyinequality} given in section \ref{EI}.

Motivated by the energy inequality Theorem \ref{energyinequality} we conclude this section by showing that the length of the set cracked the by the limiting evolution over the time interval $0\leq\tau\leq t$ is bounded. 
Recall the jump set for the deformation $u^0$ at time $\tau$ is $S_{u^0(\tau)}$ and its length is given by its one dimensional Hausdorff measure
$\mathcal{H}^1(S_{\scriptscriptstyle{u^0(\tau)}})$. The bound follows from the following theorem.

\begin{theorem}
\label{timesum}
\begin{eqnarray}
\int_0^t\left(\mathcal{G}_c \mathcal{H}^1(S_{\scriptscriptstyle{u^0(\tau)}})+\mu \int_D|\nabla u^0(\tau)|^2\,dx\right)\,d\tau<\infty, \hbox{ $0\leq t\leq T$.}
\label{timesuminequality}
\end{eqnarray}
\end{theorem}
\noindent Hence
\begin{eqnarray}
\int_0^t\mathcal{H}^1(S_{\scriptscriptstyle{u^0(\tau)}})\,d\tau<\infty, \hbox{ $0\leq t\leq T$.}
\label{timesuminequality2}
\end{eqnarray}
This shows that the total length of the set cracked by the evolution  from $t=0$ to $t=T$ is bounded. 
Theorem \ref{timesum}  is established in section \ref{SC}.

\subsection{Wave equation for the displacement}
\label{STAY}

It is shown that the limit evolution $u^0$ solves the  wave equation. The following hypothesis on the regularity of the crack set is made.
\begin{hypothesis}
\label{remarkthree}
We  suppose that the crack set given by  $S_{u^0(t)}$ is a closed set for $0\leq t\leq T$.  
\end{hypothesis}
\noindent The next hypotheses applies to the concentration set associated with unstable neighborhoods and its relation to the jump set for the limit flow.
\begin{hypothesis}
\label{remark555}
\noindent Recall from Theorem \ref{bondunstable} that the centroids of unstable neighborhoods given by Definition  \ref{UnstabeL} concentrate on the lower dimensional set $C_{0,t}$. Motivated by this observation we will assume $S_{u^0(t)}=C_{0,t}$ for $0\leq t\leq T$.
\end{hypothesis}
\noindent The next hypotheses applies to neighborhoods $\mathcal{H}_{\epsilon_k}(x)$ for which the relative displacement  is subcritical, i.e.,  $|u^{\epsilon_k}(t,y)-u^{\epsilon_k}(t,x)|<\overline{\eta}$, for $y$ in $\mathcal{H}_{\epsilon_k}(x)$. These neighborhoods will be referred to as neutrally stable.
\begin{hypothesis}
\label{remark556}
We suppose   that $\epsilon_k=\frac{1}{2^k}<\delta$ and $0\leq t\leq T$ and consider the collection of centroids $C_{\delta,t}$ associated with unstable neighborhoods. We fatten out $C_{\delta,t}$ and consider $\tilde{C}_{\delta,t}=\{x\in D:\, dist(x,C_{\delta,t})<\delta\}$. We suppose that all neighborhoods $H_{\epsilon_k}(x)$ that do not intersect the set $\tilde{C}_{\delta,t}$ are neutrally stable.
\end{hypothesis}
\noindent Passing to subsequences if necessary we apply Theorem \ref{limitflow} and  take $u^0$ to be the limit evolution of the family of  peridynamic evolutions $\{u^{\epsilon_k}\}_{k=1}^\infty$ characterized by horizons of radii $\epsilon_k=\frac{1}{2^k}$. 
\begin{theorem}{\rm \bf Wave equation.}\\
\label{waveequation}
Suppose Hypotheses \ref{remarkthree}, \ref{remark555} and \ref{remark556} hold true then the limit evolution $u^0(t,x)$ is a solution of the wave equation 
\begin{eqnarray}
\rho u^0_{tt}=2\mu\,{\rm div}(\nabla u^0)+b, \hbox{for all  $(t,x)$ on $[0,T]\times D$}.
\label{waveequationn}
\end{eqnarray}
Here the second derivative $u_{tt}^0$ is the time derivative in the sense of distributions of $u^0_t$ and ${\rm div}(\nabla u^0)$ is the divergence of the approximate gradient $\nabla u^0$ in the distributional sense.
\end{theorem}
The proof of Theorem \ref{waveequation} is given in section \ref{SC}.

\begin{remark}
\label{displacementcrack}
The sharp interface limit of the peridynamic model is given by the displacement - crack set pair $u_0(t,x)$, $S_{u^0(t)}$. 
The wave equation provides the dynamic coupling between elastic waves and the evolving fracture path inside the media.

\end{remark}

\begin{remark}
\label{remarkfinal}
We point out that the peridynamic constitutive model addressed in this work does not have an irreversibility constraint and the constitutive law \eqref{forcestate} applies at all times in the peridynamic evolution.  Because of this the crack set at each time is given by $S_{u^0(t)}$. Future work will investigate the effects of irreversibility (damage) in the peridynamic model. 
\end{remark}

\begin{remark}
\label{remarkfinalevenmoreso}
We conjecture that Hypotheses \ref{remarkthree}, \ref{remark555} and \ref{remark556} hold true. It is also pointed out that these hypotheses are only used to establish Lemma \ref{twolimitsB} which identifies the directional derivative of  approximate gradient at $x$ along the direction $e=\xi/|\xi|$ with the weak $L^2(D\times\mathcal{H}_1(0))$ limit of the difference quotients $\frac{\eta^{\epsilon_k}}{\epsilon_k|\xi|}$  restricted to pairs $(x,\xi)$ for which $|\eta^{\epsilon_k}|<\overline{\eta}$.
\end{remark}

\setcounter{equation}{0} \setcounter{theorem}{0} \setcounter{lemma}{0}\setcounter{proposition}{0}\setcounter{remark}{0}\setcounter{remark}{0}\setcounter{definition}{0}\setcounter{hypothesis}{0}

\section{Mathematical underpinnings and analysis}

From the physical perspective the convex-concave nonlinearity of the peridynamic potential  delivers the unstable constitutive law relating force to relative displacement.  On the other hand from the mathematical viewpoint this class of peridynamic potentials  share the same convex-concave structure as the function $r:\rightarrow\arctan(r^2)$ proposed by De Giorgi \cite{Gobbino1} and analyzed and generalized in the work of Gobbino \cite{Gobbino1}, \cite{Gobbino2}, and Gobbino and Mora \cite{Gobbino3} for the analysis of the Mumford Shah functional used in image processing \cite{MumfordShah}. Here we apply the methods developed in these investigations  and use them as tools for extracting the limit behavior from the peridynamic model.

In this section we provide the proofs of the theorems stated in sections two and three. The first subsection asserts the Lipschitz continuity of $\nabla PD^{\epsilon_k}(u)$ for $u$ in  $L^2_0(D)$ and applies the standard  theory of ODE to deduce existence of the peridynamic flow see, section \ref{EE}. A Gronwall inequality is used to bound the peridynamic elastic energy and kinetic energy uniformly in time see, section \ref{GKP}. We introduce $\Gamma$ -- convergence for peridynamic functions in section \ref{CC} and identify compactness conditions necessary to generate a sequence of peridynamic flows converging to a limit flow. We take limits and apply $\Gamma$ -- convergence theory to see that the  limit flows have bounded elastic energy in the sense of fracture mechanics. In section \ref{EI} we pass to the limit in the energy balance equation for peridynamic flows \eqref{BalanceEnergy} to recover an energy inequality for the limit flow. The wave equation satisfied by the limit flow is obtained on identifying the weak $L^2$ limit of the sequence $\{\nabla PD^{\epsilon_k}(u^{\epsilon_k})\}_{k=1}^\infty$  and passing to the limit in the weak formulation of \eqref{stationary} see, section \ref{SC}. We conclude with the proof of Theorem \ref{bondunstable}.

\label{UA}
\subsection{Existence of peridynamic evolution}
\label{EE}

The peridynamic equation \eqref{stationaryagain} is written as an equivalent first order system. We set $y^{\epsilon_k}=(y^{\epsilon_k}_1,y^{\epsilon_k}_2)^T$ where $y^{\epsilon_k}_1=u^{\epsilon_k}$ and $y_2^{\epsilon_k}=u_t^{\epsilon_k}$. Set $F^{\epsilon_k}(y^{\epsilon_k},t)=(F^{\epsilon_k}_1(y^{\epsilon_k},t),F^{\epsilon_k}_2(y^{\epsilon_k},t))^T$ where
\begin{eqnarray}
F^{\epsilon_k}_1(y^{\epsilon_k},t)&=&y_2^{\epsilon_k}\nonumber\\
F^{\epsilon_k}_2(y^{\epsilon_k},t)&=&\nabla PD^{\epsilon_k}(y_1^{\epsilon_k})+b(t).\nonumber
\end{eqnarray}
The initial value problem for $y^{\epsilon_k}$  given by the first order system is
\begin{eqnarray}
\frac{d}{dt} y^{\epsilon_k}=F^{\epsilon_k}(y^{\epsilon_k},t)\label{firstordersystem}
\end{eqnarray}
with initial conditions $y^{\epsilon_k}(0)=(u_0,v_0)^T$ satisfying LEFM initial conditions. In what follows we consider the more general class of initial data 
$(u_0,v_0)$ belonging to $L^2_0(D)\times L^2_0(D)$. A straight forward calculation shows that for a generic positive constant $C$ independent of $\eta$, $\xi$, and $\epsilon_k$, that
\begin{eqnarray}
\sup_{\eta}|\partial_\eta^2W^{\epsilon_k}(\eta,\epsilon_k\xi)|\leq J(|\xi|)\frac{C}{\epsilon_k^4\vert\xi|}.
\label{secondderv}
\end{eqnarray}
From this it easily follows from H\"older and Minkowski inequalities that $\nabla PD^{\epsilon_k}$ is a Lipschitz continuous map from $L^2_0(D)$ into $L^2_0(D)$ and there is a positive constant $C$ independent of $0\leq t\leq T$, such that for any pair of vectors $y=(y_1,y_2)^T$, $z=(z_1,z_2)^T$ in $L^2_0(D)\times L^2_0(D)$ 
\begin{eqnarray}
\Vert F^{\epsilon_k}(y-z,t)\Vert_{L^2(D)^2}\leq \frac{C}{\epsilon_k^2}\Vert y-z\Vert_{L^2(D)^2} \hbox{ for $0\leq t\leq T$}.
\label{lipschitz}
\end{eqnarray}
Here for any element $w=(w_1,w_2)$ of  $L^2_0(D)\times L^2_0(D)$ we have $\Vert w \Vert_{L^2(D)^2}=\Vert w_1\Vert_{L^2(D)}+\Vert w_2\Vert_{L^2(D)}$.
Since  \eqref{lipschitz} holds the standard theory of ODE in Banach space \cite{Driver} shows that there exists a unique solution to the initial value problem \eqref{firstordersystem} with $y^{\epsilon_k}$ and $\partial_t y^{\epsilon_k}$ belonging to $C([0,T]; L^2_0(D))$ and Theorem \ref{existenceuniqueness} is proved.

\subsection{Bounds on kinetic and potential energy for solutions of PD}
\label{GKP}
In this section we apply Gronwall's inequality to obtain bounds on the kinetic and elastic energy for peridynamic flows. The bounds are used to show that the solutions of the PD initial value problem are Lipschitz continuous in time.
The bounds are described in the following theorem.
\begin{theorem}
\label{Gronwall}
{\rm \bf Bounds on kinetic and potential energy for peridynamic evolution.}\\
There exists a positive constant $C$ depending only on $T$ and independent of the index $\epsilon_k$ for which
\begin{eqnarray}
\sup_{0\leq t\leq T}\left\{PD^{\epsilon_k}(u^{\epsilon_k}(t))+\frac{\rho}{2}\Vert u_t^{\epsilon_k}(t)\Vert_{L^2(D)}\right\}\leq C.
\label{boundenergy}
\end{eqnarray}
\end{theorem}

{\bf Proof.} We apply \eqref{stationaryagain} and write
\begin{eqnarray}
&&\frac{d}{dt}\left\{PD^{\epsilon_k}(u^{\epsilon_k}(t))+\frac{\rho}{2}\Vert u_t^{\epsilon_k}(t)\Vert_{L^2(D)}\right\}\nonumber\\
&&=\int_{D}(\nabla PD^{\epsilon_k}(u^{\epsilon_k}(t))+\rho u_{tt}^{\epsilon_k}(t))u_t^{\epsilon_k}(t)\,dx\nonumber\\
&&=\int_{D} u_t^{\epsilon_k}(t)b(t)\,dx\,\leq \, \frac{\rho}{2}\Vert u_t^{\epsilon_k}\Vert^2_{L^2(D)}+\frac{\rho^{-1}}{2}\Vert b(t)\Vert^2_{L^2(D)}.\label{esttime1}
\end{eqnarray}
Adding $PD^{\epsilon_k}(u^{\epsilon_k})$ to the right hand side of \eqref{esttime1} and applying Gronwall's inequality gives
\begin{eqnarray}
&& PD^{\epsilon_k}(u^{\epsilon_k}(t))+\frac{\rho}{2}\Vert u_t^{\epsilon_k}(t)\Vert_{L^2(D)}\nonumber\\
&& \leq e^t\left (PD^{\epsilon_k}(u_0)+\frac{\rho}{2}\Vert v_0\Vert_{L^2(D)}+\frac{\rho^{-1}}{2}\int_0^T\Vert b(\tau)\Vert^2_{L^2(D)}\,d\tau \right ).
\label{gineq}
\end{eqnarray}
From  \eqref{upperboundperi} of section \ref{CC} we have the upper bound
\begin{eqnarray}
PD^{\epsilon_k}(u_0)\leq LEFM (u_0,D)\hbox{ for every $\epsilon_k$, \,\,$k=1,2,\ldots$},
\label{upperbound}
\end{eqnarray}
where $LEFM(u_0,D)$ is the elastic potential energy for linear elastic fracture mechanics given by \eqref{LEFMinit}.
Theorem \ref{boundenergy} now follows from \eqref{gineq} and \eqref{upperbound}.

Theorem \ref{Gronwall} implies that PD solutions are Lipschitz continuous in time; this is stated explicitly in the following theorem. 
\begin{theorem}{\rm \bf Lipschitz continuity.}\\
\label{holder}
There exists a positive constant $K$ independent of $t_2 < t_1$ in $[0,T]$ and  index $\epsilon_k$ such that
\begin{eqnarray}
\Vert u^{\epsilon_k}(t_1)-u^{\epsilon_k}(t_2)\Vert_{L^2(D)}\leq K |t_1-t_2|.
\label{holderest}
\end{eqnarray}
\end{theorem}
{\bf Proof.} We write
\begin{eqnarray}
&&\Vert u^{\epsilon_k}(t_1)-u^{\epsilon_k}(t_2)\Vert_{L^2(D)}=\left (\int_{D}|\int_{t_1}^{t_2} u^{\epsilon_k}_\tau(\tau)\,d\tau |^2\,dx\right )^{\frac{1}{2}}\nonumber\\
&&\leq\int_{t_1}^{t_2}\Vert u_\tau^{\epsilon_k}(\tau)\Vert_{L^2(D)}\,d\tau\nonumber\\
&&\leq K|t_1-t_2|,
\label{lip}
\end{eqnarray}
where the last inequality follows from the upper bound for $\Vert u_t^{\epsilon_k}(t)\Vert_{L^2(D)}$ given by Theorem \ref{Gronwall}.

\subsection{Compactness and convergence}
\label{CC}
In this section we prove Theorem \ref{limitflow}. We start by introducing the relationship between the elastic energies $PD^{\epsilon_k}(u)$
and $LEFM(u,D)$ given by \eqref{TotalStrain} and \eqref{LEFMatt} respectively.  An application of Theorem 4.3 of \cite{Gobbino2} together with a straight forward computation using the formula for the peridynamic strain energy delivers the following inequality
\begin{eqnarray}
&&PD^{\epsilon_k}(u)\leq LEFM(u,D), \hbox{  for every $u$ in $L^2_0(D)$, and $\epsilon_k>0$}, \label{upperboundperi}
\end{eqnarray}

We now recall  the properties of $\Gamma$-convergence in order to apply them to the problem considered here. 
Consider a sequence of functions $\{F_j\}$ defined on a metric
space $\mathbb{M}$ with values in $\overline{\mathbb{R}}$ together with a function $F$ also defined on $\mathbb{M}$ with values in $\overline{\mathbb{R}}$.

\begin{definition}
\label{Gammaconvergence}
We say that $F$ is the $\Gamma$-limit of the sequence $\{F_j\}$ in $\mathbb{M}$   if the following two 
properties hold:
\begin{enumerate}
\item for every $x$ in $\mathbb{M}$ and every sequence $\{x_j\}$ converging to $x$, we have that
\begin{eqnarray}
F(x)\leq \liminf_{j\rightarrow\infty} F_j(x_j),\label{lowerbound}
\end{eqnarray}
\item for every $x$ in $\mathbb{M}$ there exists a recovery sequence $\{x_j\}$ converging to $x$, for which
\begin{eqnarray}
F(x)=\lim_{j\rightarrow\infty} F_j(x_j).\label{recovery}
\end{eqnarray}
\end{enumerate}
\end{definition}

We shall see that we can pass to the limit $\epsilon_{k}\searrow 0$ to find that the limit evolution $u^0(t,x)$ belongs to the class of Generalized SBV functions denoted by $GSBV(D)$. This class of functions has  been introduced for the study of free discontinuity problems in \cite{DeGiorgiAmbrosio} and are seen here to naturally arise in the small horizon limit of peridynamics. The space $GSBV(D)$ is composed of all measurable functions $u$ defined on $D$ whose truncations $u_k=(u\wedge k)\vee (-k)$ belong to ${SBV}(B)$
for every compact subset $B$ of $D$, see \cite{Ambrosiobook}, \cite{Buttazzo}. 
Every $u$ belonging to $GSBV(D)$ has an approximate gradient $\nabla u(x)$ for almost every $x$ in $D$ and the jump set $S_u$  is the countable union of rectifiable arcs up to a set of Hausdorff $\mathcal{H}^1$ measure zero.

For $u$ in $L^2_0(D)$ define $PD^0:\,L^2_0(D)\rightarrow [0,+\infty]$ by
\begin{equation}
PD^0(u,D)=\left\{ \begin{array}{ll}
LEFM(u,D)&\hbox{if $u$  belongs to $GSBV(D)$}\\
+\infty&\hbox{otherwise}
\end{array} \right.
\label{Gammalimit}
\end{equation}

A straight forward application of  Theorem 4.3 $(iii)$ of \cite{Gobbino2} to the sequence of peridynamic energies $\{PD^{\epsilon_k}\}$ shows that
\begin{eqnarray}
&&PD^0(u,D) \hbox{ is the $\Gamma$-limit of $\{PD^{\epsilon_k}\}$ in $L^2_0(D)$},
\label{gammaconvpd}\\
&&\lim_{k\rightarrow\infty}PD^{\epsilon_k}(u)=PD^0(u,D), \hbox{  for every $u$ in $L^2_0(D)$}.\label{pointwise}
\end{eqnarray}

Now it is shown that  the family of peridynamic flows $\{u^{\epsilon_k}\}_{k=1}^\infty$ is relatively compact in \\
$C([0,T];L^2(D))$ and
that the limit flows have bounded elastic energy in the sense of fracture mechanics.
For each $t$  in $[0,T]$ we apply Theorem \ref{Gronwall} and Hypothesis \ref{remarkone} to obtain the bound
\begin{eqnarray}
PD^{\epsilon_k}(u^{\epsilon_k}(t))+\Vert u^{\epsilon_k}(t)\Vert_{L^\infty(D)}<C
\label{cpactbnd}
\end{eqnarray}
where $C<\infty$ and is independent of $\epsilon_k$, $k=1,2,\ldots$, and $0\leq t\leq T$.
With this bound we can apply Theorem 5.1 and Remark  5.2 of \cite{Gobbino2} to assert that for each $t$ the sequence $\{u^{\epsilon_k}(t)\}_{k=1}^\infty$
is relatively compact in $L^2(D)$. 
From Theorem \ref{holder}  the sequence $\{u^{\epsilon_k}\}_{k=1}^\infty$, is seen to be uniformly equa-continuous in $t$ with respect to the $L^2(D)$ norm
and we immediately conclude from the Ascoli theorem that $\{u^{\epsilon_k}\}_{k=1}^\infty$ is relatively compact in $C([0,T];L^2(D))$.
Therefore we can pass to a  subsequence also denoted by $\{u^{\epsilon_{k}}(t)\}_{k=1}^\infty$ to assert the existence of a limit evolution $u^0(t)$ in $C([0,T];L^2(D))$  for which
\begin{eqnarray}
\lim_{k\rightarrow\infty}\left\{\sup_{t\in[0,T]}\Vert u^{\epsilon_{k}}(t)-u^0(t)\Vert_{\scriptscriptstyle{{L^2(D)}}}\right\}=0.
\label{unfconvergence}
\end{eqnarray}

Observe that since the sequence of peridynamic energies $\{PD^{\epsilon_k}\}$ $\Gamma$-converge to $PD^0$ in $L^2(D)$ we can apply the
the lower bound property \eqref{lowerbound} of $\Gamma$-convergence to conclude that the limit has bounded elastic energy in the
sense of fracture mechanics, i.e.,
\begin{eqnarray}
LEFM(u^0(t))=PD^0(u^0(t))\leq\liminf_{k\rightarrow\infty}PD^{\epsilon_{k}}(u^{\epsilon_{k}}(t))<C.
\label{GSBV}
\end{eqnarray}
This concludes the proof of Theorem \ref{limitflow}.

\subsection{Energy inequality for the limit flow}
\label{EI}
In this section we prove Theorem \ref{energyinequality}. We begin by showing that the limit evolution $u^0(t,x)$ has a weak derivative $u_t^0(t,x)$ belonging to $L^2([0,T]\times D)$.  This is summarized in the following theorem.
\begin{theorem}
\label{weaktimederiviative}
On passage to subsequences as necessary the sequence $u_t^{\epsilon_k}$ weakly converges in $L^2([0,T]\times D)$ to $u^0_t$ where
\begin{eqnarray}
-\int_0^T\int_D\partial_t\psi u^0\, dxdt=\int_0^T\int_D\psi u^0_t\, dxdt,
\label{weakl2time}
\end{eqnarray}
for all compactly supported smooth test functions $\psi$ on $[0,T]\times D$.
\end{theorem}

{\bf Proof.} The bound on the kinetic energy given in Theorem \ref{Gronwall} 
implies
\begin{eqnarray}
\sup_{\epsilon_k>0}\left(\sup_{0\leq t\leq T}\Vert u^{\epsilon_k}_t\Vert_{L^2(D)}\right)< \infty.
\label{bddd}
\end{eqnarray}
Therefore the sequence $u^{\epsilon_k}_t$ is bounded in $L^2([0,T]\times D)$ and passing to
a subsequence if necessary we conclude that there is a limit function
$\tilde{u}^0$ for which $u_t^{\epsilon_k}\rightharpoonup\tilde{u}^0$ weakly in $L^2([0,T]\times D)$. Observe also that the uniform convergence \eqref{unfconvergence} implies that $u^{\epsilon_k}\rightarrow u^0$ in $L^2([0,T]\times D)$. 
On writing the identity 
\begin{eqnarray}
-\int_0^T\int_D\partial_t\psi u^{\epsilon_k}\, dxdt=\int_0^T\int_D\psi u^{\epsilon_k}_t\, dxdt.
\label{weakidentity}
\end{eqnarray}
applying our observations and passing to the limit it is seen that
$\tilde{u}^0=u_t^0$ and the theorem follows.

To establish Theorem \ref{energyinequality} we require the following inequality.
\begin{lemma}
For every $t$ in $[0,T]$ we have
\label{weakinequality}
\begin{eqnarray}
\Vert u^0_t(t)\Vert_{L^2(D)}\leq \liminf_{\epsilon_k\rightarrow 0}\Vert u^{\epsilon_k}_t(t)\Vert_{L^2(D)}.
\label{limitweakineq}
\end{eqnarray}
\end{lemma}
{\bf Proof.}
For every non-negative bounded measurable  function of time $\psi(t)$ defined on $[0,T]$ we have the inequality
\begin{eqnarray}
\int_0^t\psi \Vert u^{\epsilon_k}_t-u^0_t\Vert_{L^2(D)}^2\,dt\geq 0
\label{positive}
\end{eqnarray}
and together with the weak convergence given in Theorem \ref{weaktimederiviative}  one easily sees that
\begin{eqnarray} 
\liminf_{\epsilon_k\rightarrow 0}\int_0^T\psi\Vert u^{\epsilon_k}_t\Vert_{L^2(D)}^2\,dt-\int_0^T\psi\Vert u^0_t\Vert_{L^2(D)}^2\,dt\geq 0.
\label{diff}
\end{eqnarray}

\noindent Applying  \eqref{bddd} and invoking the Lebesgue dominated convergence theorem we conclude
\begin{eqnarray}
\liminf_{\epsilon_k\rightarrow 0}\int_0^T\psi\Vert u^{\epsilon_k}_t\Vert_{L^2(D)}^2\,dt=\int_0^T\psi\liminf_{\epsilon_k\rightarrow 0}\Vert u^{\epsilon_k}_t\Vert_{L^2(D)}^2\,dt
\label{equalitybalance}
\end{eqnarray}
to recover the inequality given by
\begin{eqnarray}
\int_0^T\psi\left(\liminf_{\epsilon_k\rightarrow 0}\Vert u^{\epsilon_k}_t\Vert_{L^2(D)}^2-\Vert u^0_t\Vert_{L^2(D)}^2\right)\,dt\geq 0.
\label{diffinal}
\end{eqnarray}
The lemma follows noting that \eqref{diffinal} holds for every non-negative test function $\psi$.

Theorem \ref{energyinequality} now follows immediately on taking the $\epsilon_k\rightarrow 0$ limit in the peridynamic energy balance equation \eqref{BalanceEnergy} of Theorem \ref{Ebalance} and applying \eqref{pointwise}, \eqref{unfconvergence}, \eqref{GSBV}, and \eqref{limitweakineq} of Lemma \ref{weakinequality}.

\subsection{Stationarity conditions for the limit flow}
\label{SC}

In this section we prove Theorems \ref{timesum} and \ref{waveequation}. In the first subsection we give the proof of Theorem \ref{timesum}.
In the second subsection we provide the proof of Theorem \ref{waveequation} using Theorem \ref{convgofelastic}. In the last subsection we prove Theorem \ref{convgofelastic} .

\subsubsection{Proof of Theorem \ref{timesum}}

We consider the integral 
\begin{eqnarray}
\label{integralboundfirst}
\int_0^t\left(PD^{\epsilon_k}(u^{\epsilon_k}(\tau)+\Vert u_t^{\epsilon_k}(\tau)\Vert_{L^2(D)}^2\right)\,d\tau,
\end{eqnarray}
and apply the energy bound \eqref{boundenergy} to obtain the inequality
\begin{eqnarray}
\label{integralboundtime}
\int_0^t\left(PD^{\epsilon_k}(u^{\epsilon_k}(\tau)+\Vert u_t^{\epsilon_k}(\tau)\Vert_{L^2(D)}^2\right)\,d\tau<Ct.
\end{eqnarray}
Since $PD^{\epsilon_k}(u^{\epsilon_k}(t))$ and $\Vert u_t^{\epsilon_k}(t)\Vert_{L^2(D)}^2$ are non-negative we
apply Fatou's Lemma to see that
\begin{eqnarray}
\label{integralboundtimefatou}
\int_0^t\liminf_{\epsilon_k\rightarrow 0}\left(PD^{\epsilon_k}(u^{\epsilon_k}(\tau)+\Vert u_t^{\epsilon_k}(\tau)\Vert_{L^2(D)}^2\right)\,d\tau<Ct.
\end{eqnarray}
Applying \eqref{GSBV} and \eqref{limitweakineq} delivers the upper bound
\begin{eqnarray}
\label{upperboundsurface}
\int_0^t  LEFM(u^0(\tau),D) \,d\tau+\int_0^t\Vert u^0_t(\tau)\Vert_{L^2(D)}^2\,d\tau<Ct.
\end{eqnarray}
Here we have used the fact that $PD^{\epsilon_k}(u^0(t))$ is continuous in $t$ and the pointwise convergence $PD^{\epsilon_k}(u^0(t))\rightarrow LEFM(u^0(t),D)$ to assert the integrability of $LEFM(u^0(t),D)$ with respect to $t$.
Theorem \ref{timesum} now follows from \eqref{upperboundsurface}.

\subsubsection{Proof of Theorem \ref{waveequation}}

We introduce the following integration by parts identity that holds for any pair of functions $u$, $v$ belonging to $L_0^2(D)$ with either $u$ or $v$ having compact support inside $D$ given by
\begin{eqnarray}
\int_{D}\int_{\mathcal{H}_1(0)} D_{-e}^{{\epsilon_k}|\xi|} u v\,d\xi\,dx=\int_{D}\int_{\mathcal{H}_1(0)} u D_e^{{\epsilon_k}|\xi|} v\,d\xi\,dx.
\label{intbyparts}
\end{eqnarray}
Note further if $v$ is infinitely differentiable and has compact support in $D$ then
\begin{eqnarray}
&&\lim_{\epsilon_k\rightarrow 0}D_e^{{\epsilon_k}|\xi|} v=\nabla v\cdot e\label{grad}
\end{eqnarray}
where the convergence is uniform in $D$. Here $e$ is the unit vector $e=\xi/|\xi|$.

Taking the first variation of the action integral \eqref{Action} gives the Euler equation in weak form
\begin{eqnarray}
\rho\int_0^T\int_{D} u_t^{\epsilon_k}\delta_t\,dx \,dt-\int_0^T\int_{D}\nabla PD^{\epsilon_k}(u^{\epsilon_k})\delta\,dx\,dt+\int_0^T\int_{D} b\delta\,dx\,dt=0\label{stationtx}
\end{eqnarray}
where the test function $\delta=\delta(x,t)=\psi(t)\phi(x)$ is smooth and has compact support in $[0,T]\times D$. Integrating by parts in the second term of \eqref{stationtx} using \eqref{intbyparts} gives
\begin{eqnarray}
&&\rho\int_0^T\int_{D} u_t^{\epsilon_k}\delta_t\,dx \,dt\nonumber\\
&&-\int_0^T\int_{D}\int_{\mathcal{H}_1(0)}|\xi|J(|\xi|)f'\left(\frac{|\eta^{\epsilon_k}|^2}{\epsilon_k|\xi|}\right)\frac{2\eta^{\epsilon_k}}{\epsilon|\xi|}D_{e}^{\epsilon_k}\delta\,d\xi\,dx\,dt+\int_0^T\int_{D} b\delta\,dx\,dt=0.\label{stationtxweaker}
\end{eqnarray}
Where $\eta^{\epsilon_k}=u^{\epsilon_k}(x+\xi)-u^{\epsilon_k}(x)$ and observe that $\eta^{\epsilon_k}/(\epsilon_k |\xi|)=D_e^{{\epsilon_k}|\xi|} u^{\epsilon_k}$. Next we make the change of function and write $F_s (r)=\frac{1}{s}f(sr^2)$ and on setting $s={\epsilon_k}|\xi|$ and $r=D_e^{{\epsilon_k}|\xi|} u^{\epsilon_k}$ we transform \eqref{stationtxweaker} into
\begin{eqnarray}
&&\rho\int_0^T\int_{D} u_t^{\epsilon_k}\delta_t\,dx \,dt\nonumber\\
&&-\int_0^T\int_{D}\int_{\mathcal{H}_1(0)}|\xi|J(|\xi|)F_{\epsilon_k|\xi|}'(D_{e}^{{\epsilon_k}|\xi|}u^{\epsilon_k})D_{e}^{\epsilon_k}\delta\,d\xi\,dx\,dt+\int_0^T\int_{D} b\delta\,dx\,dt=0,\label{stationtxweakerlimitform}
\end{eqnarray}
where
\begin{eqnarray}
F_{\epsilon_k|\xi|}'(D_{e}^{\epsilon_k |\xi|}u^{\epsilon_k})=f'\left(\epsilon_k |\xi||D_{e}^{\epsilon_k |\xi|}u^{\epsilon_k}|^2\right)2D_{e}^{\epsilon_k |\xi|}u^{\epsilon_k}.
\label{deriv}
\end{eqnarray}

For future reference observe that $F_s(r)$ is convex-concave in $r$ with inflection point $\overline{r}_s=\overline{r}/\sqrt{s}$ where $\overline{r}$ is the inflection point of $f(r^2)=F_1(r)$. One also has the estimates
\begin{eqnarray}
&&F_s(r)\geq\frac{1}{s}F_1(\overline{r})\hbox{  for $r\geq\overline{r}_s$, and }\label{lowerestforF}\\
&&\sup_{0\leq r<\infty}|F'_s(r)|\leq\frac{2f'(\overline{r}^2){\overline{r}}}{\sqrt{s}},\label{boundderiv}
\label{estforFprime}
\end{eqnarray}
We send $\epsilon_k\rightarrow 0$ in \eqref{stationtxweakerlimitform} applying the weak convergence Theorem \ref{weaktimederiviative} to the first term to obtain
\begin{eqnarray}
&&\rho\int_0^T\int_{D} u_t^{0}\delta_t\,dx \,dt-\lim_{\epsilon_k\rightarrow 0}\left(\int_0^T\int_{D}\int_{\mathcal{H}_1(0)}|\xi|J(|\xi|)F_{\epsilon_k|\xi|}'(D_{e}^{{\epsilon_k}|\xi|}u^{\epsilon_k})D_{e}^{\epsilon_k}\delta\,d\xi\,dx\,dt\right)\nonumber\\
&&+\int_0^T\int_{D} b\delta\,dx\,dt=0,\label{stationtxweakerlimit}
\end{eqnarray}
Theorem \ref{waveequation} follows once we identify the limit of the second term in \eqref{stationtxweakerlimit} for smooth test functions $\phi(x)$ with support contained in $D$.
We state the following convergence theorem.
\begin{theorem}
\label{convgofelastic}
Given any infinitely differentiable test function $\phi$ with compact support in $D$ then
\begin{eqnarray}
\lim_{\epsilon_k\rightarrow 0}\int_{D}\int_{\mathcal{H}_1(0)}|\xi|J(|\xi|)F_{\epsilon_k|\xi|}'(D_{e}^{{\epsilon_k}|\xi|}u^{\epsilon_k})D_{e}^{\epsilon_k}\phi\,d\xi\,dx=2\mu\int_{D}\nabla\phi\cdot\nabla u^0\,dx,
\label{limitdelta}
\end{eqnarray}
where $\mu=\pi f'(0)\int_0^1 r^2J(r)dr$.
\end{theorem}
\noindent Theorem \ref{convgofelastic} is proved in section \ref{prtheorem44}. 
The sequence of integrals on the left hand side of \eqref{limitdelta} are uniformly bounded in time, i.e.,
\begin{eqnarray}
\sup_{\epsilon_k>0}\left\{\sup_{0\leq t\leq T}\left\vert\int_{D}\int_{\mathcal{H}_1(0)}|\xi|J(|\xi|)F_{\epsilon_k|\xi|}'(D_{e}^{{\epsilon_k}|\xi|}u^{\epsilon_k})D_{e}^{\epsilon_k}\phi\,d\xi\,dx\right\vert\right\}<\infty,
\label{uniboundt}
\end{eqnarray}
this is demonstrated in \eqref{Fprimesecond} of Lemma \ref{estimates}  in section \ref{prtheorem44}. Applying the Lebesgue bounded convergence theorem together with Theorem \ref{convgofelastic}
with $\delta(t,x)=\psi(t)\phi(x)$ delivers the desired result
\begin{eqnarray}
&&\lim_{\epsilon_k\rightarrow 0}\left(\int_0^T\int_{D}\int_{\mathcal{H}_1(0)}|\xi|J(|\xi|)F_{\epsilon_k|\xi|}'(D_{e}^{{\epsilon_k}|\xi|}u^{\epsilon_k})\psi D_{e}^{\epsilon_k}\phi\,d\xi\,dx\,dt\right)\nonumber\\
&&=2\mu\int_0^T\int_{D}\psi\nabla\phi\cdot\nabla u^0\,dx\,dt,
\label{limitidentity}
\end{eqnarray}
and we recover the identity
\begin{eqnarray}
&&\rho\int_0^T\int_{D} u_t^{0}(t,x)\psi_t(t)\phi(x)\,dx \,dt-2\mu\int_0^T\int_{D}\psi(t)\nabla\phi(x)\cdot\nabla u^0(t,x)\,dx\,dt\nonumber\\
&&+\int_0^T\int_{D} b(t,x)\psi(t)\phi(x)\,dx\,dt=0
\label{finalweakidentity}
\end{eqnarray}
from which Theorem \ref{waveequation} follows.

\subsubsection{ Proof of Theorem \ref{convgofelastic}}
\label{prtheorem44}

We decompose the difference $D_e^{\epsilon_k |\xi|}u^{\epsilon_k}$ as 
\begin{eqnarray}
D_e^{\epsilon_k |\xi|}u^{\epsilon_k}=D_e^{\epsilon_k |\xi|,-}u^{\epsilon_k}+D_e^{\epsilon_k |\xi|,+}u^{\epsilon_k}
\label{decompose}
\end{eqnarray}
where
\begin{equation}
D_e^{\epsilon_k |\xi|,-}u^{\epsilon_k}=\left\{\begin{array}{ll}D_e^{\epsilon_k |\xi|}u^{\epsilon_k},&\hbox{if  $|D_e^{\epsilon_k |\xi|}u^{\epsilon_k}|<\frac{\overline{r}}{\sqrt{\epsilon_k|\xi|}}$}\\
0,& \hbox{otherwise}
\end{array}\right.
\label{decomposedetails}
\end{equation}
where $\overline{r}$ is the inflection point for the function $F_1(r)=f(r^2)$. Here $D_e^{\epsilon_k |\xi|,+}u^{\epsilon_k}$ is defined so that \eqref{decompose} holds.
We prove Theorem \ref{convgofelastic} by using the following two identities
described in the Lemmas below.
\begin{lemma}
\label{twolimitsA}
For any $\psi$ in $C^\infty_0(D)$ 
\begin{eqnarray}
&&\lim_{\epsilon_k\rightarrow 0} \int_{D}\int_{\mathcal{H}_1(0)}|\xi|J(|\xi|)F_{\epsilon_k|\xi|}'(D_{e}^{{\epsilon_k}|\xi|}u^{\epsilon_k})D_{e}^{\epsilon_k|\xi|}\psi\,d\xi\,dx\nonumber\\
&&-2\lim_{\epsilon_k\rightarrow 0}\int_{D}\int_{\mathcal{H}_1(0)}|\xi|J(|\xi|)f'(0)D_e^{\epsilon_k |\xi|,-}u^{\epsilon_k}D_e^{\epsilon_k |\xi|}\psi\,d\xi\,dx=0.
\label{functiotolimit}
\end{eqnarray}
\end{lemma}
\begin{lemma}
\label{twolimitsB}
Assume that Hypotheses \ref{remarkthree}, \ref{remark555} and \ref{remark556} hold true and
define the weighted Lebesgue measure $\nu$ by $\nu(S)=\int_S|\xi|J(|\xi|)d\xi\,dx$ for any Lebesgue measurable set $S\subset D\times\mathcal{H}_1(0)$.
Passing to subsequences as necessary  $\{D_e^{\epsilon_k |\xi|,-}u^{\epsilon_k}\}_{k=1}^\infty$ converges weakly in $L^2(D\times\mathcal{H}_1(0);\nu)$ to $e\cdot \nabla u^0$ where $e=\xi/|\xi|$, i.e., 
\begin{eqnarray}
&&\lim_{\epsilon_k\rightarrow 0}\int_{D}\int_{\mathcal{H}_1(0)}|\xi|J(|\xi|)D_e^{\epsilon_k |\xi|,-}u^{\epsilon_k} \phi\,d\xi\,dx\nonumber\\
&&=\int_{D}\int_{\mathcal{H}_1(0)}|\xi|J(|\xi|)e\cdot\nabla u^0\phi\,d\xi\,dx,
\label{functiotofunction}
\end{eqnarray}
for any test function $\phi(x,\xi)$ in  $L^2(D\times\mathcal{H}_1(0);\nu)$. 
\end{lemma}

We now apply the Lemmas. Observing that $D_e^{\epsilon_k |\xi|}\psi$ converges strongly in $L^2(D\times\mathcal{H}_1(0):\nu)$ to $e\cdot\nabla\psi$ for test functions $\psi$ in $C^\infty_0(D)$ and from the weak $L^2(D\times\mathcal{H}_1(0):\nu)$ convergence of $D_e^{\epsilon_k |\xi|,-}u^{\epsilon_k}$ we deduce that
\begin{eqnarray}
&&\lim_{\epsilon_k\rightarrow 0}\int_{D}\int_{\mathcal{H}_1(0)}|\xi|J(|\xi|)f'(0)D_e^{\epsilon_k |\xi|,-}u^{\epsilon_k}D_e^{\epsilon_k |\xi|}\psi\,d\xi\,dx\nonumber\\
&&=\int_{D}\int_{\mathcal{H}_1(0)}|\xi|J(|\xi|)f'(0)\,(e\cdot\nabla u^0)(e\cdot\nabla\psi)\,d\xi\,dx\nonumber\\
&&=f'(0)\int_{\mathcal{H}_1(0)}|\xi|J(|\xi|)\,e_i e_j\,d\xi\int_{D}\partial_{x_i}u^0\partial_{x_i}\psi\,dx.\label{limitproduct}
\end{eqnarray}
A calculation shows that
\begin{eqnarray}
f'(0)\int_{\mathcal{H}_1(0)}|\xi|J(|\xi|)\,e_i e_j\,d\xi=\mu\delta_{ij}\label{shear}
\end{eqnarray}
where $\mu$ is given by \eqref{calibrate}.
Theorem \ref{convgofelastic} now follows immediately from \eqref{limitproduct} and \eqref{functiotolimit}. 

To establish Lemmas \ref{twolimitsA} and \ref{twolimitsB}  we develop the following estimates for the sequences $D_e^{\epsilon_k |\xi|,-}u^{\epsilon_k}$ and $D_e^{\epsilon_k |\xi|,+}u^{\epsilon_k}$.  We define the set $K^{+,\epsilon_k}$ by
\begin{eqnarray}
K^{+,\epsilon_k}=\{(x,\xi)\in D\times\mathcal{H}_1(0)\,: D_e^{\epsilon_k |\xi|,+}u^{\epsilon_k}\not=0\}.
\label{suppset}
\end{eqnarray}

We have the following string of estimates.
\begin{lemma}
We introduce the generic positive constant $0<C<\infty$ independent of $0<\epsilon_k<1$ and $0\leq t\leq T$ and state the following inequalities that hold for all $0<\epsilon_k<1$ and $0\leq t\leq T$ and for $C^\infty(D)$ test functions $\phi$  with compact support on $D$.
\label{estimates}
\begin{eqnarray}
&&\int_{K^{+,\epsilon_k}} |\xi|J(|\xi|)\,d\xi\,dx<C\epsilon_k,\label{suppsetupper}\\
&&\left\vert\int_{D\times\mathcal{H}_1(0)}|\xi| J(|\xi|)F_{\epsilon_k|\xi|}'(D_{e}^{{\epsilon_k}|\xi|,+}u^{\epsilon_k})D_{e}^{\epsilon_k}\phi\,d\xi\,dx\right\vert<C\sqrt{\epsilon_k}\Vert\nabla \phi\Vert_{L^\infty(D)},\label{Fprimefirst}\\
&&\int_{D\times\mathcal{H}_1(0)}|\xi| J(|\xi|)|D_{e}^{{\epsilon_k}|\xi|,-}u^{\epsilon_k}|^2\,d\xi\,dx<C,\label{L2bound}\\
&&\int_{D\times\mathcal{H}_1(0)}|\xi| J(|\xi|)|D_{e}^{{\epsilon_k}|\xi|}u^{\epsilon_k}|\,d\xi\,dx<C,\hbox{  and}\label{L1bound}\\
&&\left\vert\int_{D\times\mathcal{H}_1(0)}|\xi| J(|\xi|)F_{\epsilon_k|\xi|}'(D_{e}^{{\epsilon_k}|\xi|}u^{\epsilon_k}) D_{e}^{\epsilon_k}\phi\,d\xi\,dx\right\vert<C\Vert\nabla \phi\Vert_{L^\infty(D)}.\label{Fprimesecond}
\end{eqnarray}
\end{lemma}
{\bf Proof.}
For $(x,\xi)\in K^{+,\epsilon_k}$ we apply \eqref{lowerestforF} to get
\begin{eqnarray}
J(|\xi|)\frac{1}{\epsilon_k}F_{1}(\overline{r})=|\xi|J(|\xi|)\frac{1}{\epsilon_k|\xi|}F_{1}(\overline{r})
\leq|\xi|J(|\xi|)F_{\epsilon_k|\xi|}(D_{e}^{{\epsilon_k}|\xi|}u^{\epsilon_k})\label{first}
\end{eqnarray}
and in addition since $|\xi|\leq 1$ we have
\begin{eqnarray}
&&\frac{1}{\epsilon_k}F_{1}(\overline{r}) \int_{K^{+,\epsilon_k}} |\xi|J(|\xi|)\,d\xi\,dx\leq\frac{1}{\epsilon_k}F_{1}(\overline{r}) \int_{K^{+,\epsilon_k}} J(|\xi|)\,d\xi\,dx\nonumber\\
&&\leq\int_{K^{+,\epsilon_k}} |\xi|J(|\xi|)F_{\epsilon_k|\xi|}(D_{e}^{{\epsilon_k}|\xi|}u^{\epsilon_k})\,d\xi\,dx\leq\sup_{t\in [0,T]}\sup_{\epsilon_k}PD^{\epsilon_k}(u^{\epsilon_k})
\label{firstb}
\end{eqnarray}
where Theorem \ref{Gronwall} implies that the right most element of  the sequence of inequalities is bounded and \eqref{suppsetupper} follows noting that the inequality \eqref{firstb} is equivalent to \eqref{suppsetupper}. More generally since $|\xi|\leq 1$ we may argue as above to conclude that 
\begin{eqnarray}
\int_{K^{+,\epsilon_k}} |\xi|^pJ(|\xi|)\,d\xi\,dx<C\epsilon_k.
\label{power}
\end{eqnarray}
for $0\leq p$.
We apply \eqref{estforFprime} and \eqref{power} to find
\begin{eqnarray}
&&\left\vert\int_{D\times\mathcal{H}_1(0)}|\xi| J(|\xi|)F_{\epsilon_k|\xi|}'(D_{e}^{{\epsilon_k}|\xi|,+}u^{\epsilon_k})D_{e}^{\epsilon_k}\phi\,d\xi\,dx\right\vert\nonumber\\
&&\leq\frac{2f'(\overline{r}^2)\overline{r}}{\sqrt{\epsilon_k}}\int_{K^{+,\epsilon_k}}\sqrt{|\xi|}J(|\xi|)\,d\xi\,dx
\leq\sqrt{\epsilon_k}C,
\label{second}
\end{eqnarray}
and \eqref{Fprimefirst} follows.

A basic calculation shows there exists a positive constant independent of $r$ and $s$ for which
\begin{eqnarray}
r^2\leq C F_s(r), \hbox{  for $r<\frac{\overline{r}}{\sqrt{s}}$},
\label{rsquared}
\end{eqnarray}
so
\begin{eqnarray}
|D_{e}^{\epsilon_k |\xi|}u^{\epsilon_k}|^2\leq C F_{\epsilon_k|\xi|}(D_{e}^{{\epsilon_k}|\xi|}u^{\epsilon_k}), \hbox{  for $|D_{e}^{\epsilon_k|\xi|}u^{\epsilon_k}|<\frac{\overline{r}}{\sqrt{\epsilon_k |\xi|}}$},
\label{rsquared}
\end{eqnarray}
and
\begin{eqnarray}
&&\int_{D\times\mathcal{H}_1(0)}|\xi| J(|\xi|)|D_{e}^{{\epsilon_k}|\xi|,-}u^{\epsilon_k}|^2\,d\xi\,dx
=\int_{D\times\mathcal{H}_1(0)\setminus K^{+,\epsilon_k}}|\xi| J(|\xi|)|D_{e}^{{\epsilon_k}|\xi|}u^{\epsilon_k}|^2\,d\xi\,dx\nonumber\\
&&\leq C\int_{D\times\mathcal{H}_1(0)\setminus K^{+,\epsilon_k}}|\xi| J(|\xi|)F_{\epsilon_k |\xi|}(D_{e}^{\epsilon_k |\xi|}u^{\epsilon_k})\,d\xi\,dx\leq C\sup_{t\in [0,T]}\sup_{\epsilon_k}PD^{\epsilon_k}(u^{\epsilon_k})
\label{third}
\end{eqnarray}
where Theorem \ref{Gronwall} implies that the right most element of  the sequence of inequalities is bounded and \eqref{L2bound} follows.

To establish \eqref{L1bound} we apply H\"olders inequality to find that
\begin{eqnarray}
&&\int_{D\times\mathcal{H}_1(0)}|\xi| J(|\xi|)|D_{e}^{{\epsilon_k}|\xi|}u^{\epsilon_k}|\,d\xi\,dx\nonumber\\
&&=\int_{K^{+,\epsilon_k}}|\xi| J(|\xi|)|D_{e}^{{\epsilon_k}|\xi|}u^{\epsilon_k}|\,d\xi\,dx+\int_{D\times\mathcal{H}_1(0)\setminus K^{+,\epsilon_k}}|\xi| J(|\xi|)|D_{e}^{{\epsilon_k}|\xi|}u^{\epsilon_k}|\,d\xi\,dx\nonumber\\
&&\leq \frac{2\Vert u^{\epsilon_k}\Vert_{L^\infty(D)}}{\epsilon_k}\int_{K^{+,\epsilon_k}}J(|\xi|)\,d\xi\,dx+\nonumber\\
&&+\nu(D\times\mathcal{H}_1(0))^{\frac{1}{2}}\left (\int_{D\times\mathcal{H}_1(0)}|\xi| J(|\xi|)|D_{e}^{{\epsilon_k}|\xi|,-}u^{\epsilon_k}|^2\,d\xi\,dx\right)^{\frac{1}{2}},
\label{twoterms}
\end{eqnarray}
and \eqref{L1bound} follows from \eqref{power} and \eqref{L2bound}, and \eqref{max}.

We establish \eqref{Fprimesecond}. This bound follows from the basic features of the potential function $f$. We will recall for subsequent use that $f$ is smooth positive, concave and $f'$ is a decreasing function with respect to its argument. So for $A$ fixed and $0\leq h\leq A^2\overline{r}^2$ we have
\begin{eqnarray}
|f'(h)-f'(0)|\leq |f'(A^2\overline{r}^2)- f'(0)|<2|f'(0)|^2.
\label{ffact}
\end{eqnarray}
The bound \eqref{Fprimesecond} is now shown to be a consequence of the following upper bound valid for the parameter $0<A<1$ given by
\begin{eqnarray}
&&\int_{D\times\mathcal{H}_1(0)}|\xi|J(|\xi|)|f'(\epsilon_k |\xi||D_{e}^{\epsilon_k |\xi|,-}u^{\epsilon_k}|^2)-f'(0)|^2\, d\xi\,dx\nonumber\\
&&\leq \nu(D\times\mathcal{H}_1(0))\times|f'(A^2\overline{r}^2)-f'(0)|^2+C\epsilon_k\frac{4|f'(0)|^2}{A^2}.
\label{usefulbound}
\end{eqnarray}
We postpone the proof of \eqref{usefulbound} until after it is used to establish  \eqref{Fprimesecond}. Set $h_{\epsilon_k}=D_{e}^{\epsilon_k |\xi|,-}u^{\epsilon_k}$ to note 
\begin{eqnarray}
F_{\epsilon_k |\xi|}'(h_{\epsilon_k})-2f'(0)h_{\epsilon_k}=(f'(\epsilon_k |\xi| h^2_{\epsilon_k})-f'(0))2h_{\epsilon_k}.
\label{diffeq}
\end{eqnarray}
Applying H\"olders inequality, \eqref{Fprimefirst}, \eqref{L2bound}, \eqref{diffeq}, and \eqref{usefulbound} gives
\begin{eqnarray}
&&\left\vert\int_{D\times\mathcal{H}_1(0)}|\xi| J(|\xi|)F_{\epsilon_k|\xi|}'(D_{e}^{{\epsilon_k}|\xi|}u^{\epsilon_k}) D_{e}^{\epsilon_k}\phi\,d\xi\,dx
\right\vert\nonumber\\
&&\leq\left\vert\int_{D\times\mathcal{H}_1(0)}|\xi| J(|\xi|)F_{\epsilon_k|\xi|}'(D_{e}^{{\epsilon_k,+}|\xi|}u^{\epsilon_k}) D_{e}^{\epsilon_k}\phi\,d\xi\,dx
\right\vert\nonumber\\
&&+\left\vert\int_{D\times\mathcal{H}_1(0)}|\xi| J(|\xi|)F_{\epsilon_k|\xi|}'(D_{e}^{{\epsilon_k,-}|\xi|}u^{\epsilon_k}) D_{e}^{\epsilon_k}\phi\,d\xi\,dx
\right\vert\nonumber\\
&&\leq C\sqrt{\epsilon_k}\Vert\nabla \phi\Vert_{L^\infty(D)}+2\int_{D\times \mathcal{H}_1(0)}|\xi|J(|\xi|)f'(0)D_e^{\epsilon_k |\xi|,-}u^{\epsilon_k}D_e^{\epsilon_k |\xi|}\psi\,d\xi\,dx\nonumber\\
&&+\int_{D\times\mathcal{H}_1(0)}|\xi|J(|\xi|)\left (F_{\epsilon_k|\xi|}'(D_{e}^{{\epsilon_k}|\xi|}u^{\epsilon_k}) -2f'(0)D_e^{\epsilon_k |\xi|,-}u^{\epsilon_k}\right )D_e^{\epsilon_k |\xi|}\psi\,d\xi\,dx\nonumber\\
&&\leq 2C\left(f'(0)+\sqrt{\epsilon_k}+\nu(D\times\mathcal{H}_1(0))\times|f'(A^2\overline{r}^2)-f(0)|^2+\epsilon_k\frac{4|f'(0)|^2}{A^2}\right )\Vert\nabla \phi\Vert_{L^\infty(D)}.\nonumber\\
\label{five}
\end{eqnarray}
and \eqref{Fprimesecond} follows.

We establish the inequality \eqref{usefulbound}. Set $h_{\epsilon_k}=D_{e}^{\epsilon_k |\xi|,-}u^{\epsilon_k}$ and for $0<A<1$ introduce the set
\begin{eqnarray}
K^{+,\epsilon_k}_A=\{(x,\xi)\in D\times\mathcal{H}_1(0)\,: A^2\overline{r}^2\leq \epsilon_k|\xi||h_{\epsilon_k}|^2\}.
\label{suppAset}
\end{eqnarray}
To summarize  $(x,\xi)\in K^{+,\epsilon_k}_A$ implies $A^2\overline{r}^2\leq\epsilon_k|\xi||h_{\epsilon_k}|^2\leq\overline{r}^2$ and $(x,\xi)\not\in K^{+,\epsilon_k}_A$ implies $\epsilon_k|\xi||h_{\epsilon_k}|^2<A^2\overline{r}^2$ and $|f'(\epsilon_k|\xi||h_{\epsilon_k}|^2)-f'(0)|\leq|f'(A^2\overline{r}^2)-f'(0)|$. Inequality \eqref{L2bound} implies
\begin{eqnarray}
&&C>\int_{K^{+,\epsilon_k}_A} |\xi|J(|\xi|) h_{\epsilon_k}^2\,d\xi\,dx\geq\frac{A^2\overline{r}^2}{\epsilon_k}\int_{K^{+,\epsilon_k}_A} J(|\xi|) \,d\xi\,dx\nonumber\\
&&\geq\frac{A^2\overline{r}^2}{\epsilon_k}\int_{K^{+,\epsilon_k}_A} |\xi|J(|\xi|) \,d\xi\,dx,
\label{chebyA}
\end{eqnarray}
the last inequality follows since $1\geq|\xi|>0$. Hence
\begin{eqnarray}
\int_{K^{+,\epsilon_k}_A} |\xi|J(|\xi|) \,d\xi\,dx\leq C\frac{\epsilon_k}{A^2\overline{r}^2},
\label{chebyAUpper}
\end{eqnarray}
and it follows that
\begin{eqnarray}
&&\int_{K^{+,\epsilon_k}_A} |\xi|J(|\xi|)|f'(\epsilon_k|\xi|h_{\epsilon_k}|^2-f'(0)|^2 \,d\xi\,dx\nonumber\\
&&\leq 4|f'(0)|^2\int_{K^{+,\epsilon_k}_A} |\xi|J(|\xi|) \,d\xi\,dx\leq C\epsilon_k\frac{4|f'(0)|^2}{A^2\overline{r}^2}.
\label{kepsplus}
\end{eqnarray}
Collecting observations gives
\begin{eqnarray}
&&\int_{D\times\mathcal{H}_1(0)\setminus K^{+,\epsilon_k}_A}|\xi|J(|\xi|)|f'(\epsilon_k |\xi||D_{e}^{\epsilon_k |\xi|,-}u^{\epsilon_k}|^2)-f'(0)|^2\, d\xi\,dx\nonumber\\
&&\leq \nu(D\times\mathcal{H}_1(0))\times |f'(A^2\overline{r}^2)-f'(0)|^2,
\label{2ndbd}
\end{eqnarray}
and \eqref{usefulbound} follows.

We now prove Lemma \ref{twolimitsA}. Write
\begin{eqnarray}
F_{\epsilon_k|\xi|}'(D_{e}^{\epsilon_k |\xi|}u^{\epsilon_k})=F_{\epsilon_k|\xi|}'(D_{e}^{\epsilon_k |\xi|,+}u^{\epsilon_k})
+F_{\epsilon_k|\xi|}'(D_{e}^{\epsilon_k |\xi|,-}u^{\epsilon_k}),
\label{plusminus}
\end{eqnarray}
and from \eqref{Fprimefirst} it follows that
\begin{eqnarray}
&&\lim_{\epsilon_k\rightarrow 0} \int_{D}\int_{\mathcal{H}_1(0)}|\xi|J(|\xi|)F_{\epsilon_k|\xi|}'(D_{e}^{{\epsilon_k}|\xi|}u^{\epsilon_k})D_{e}^{\epsilon_k|\xi|}\psi\,d\xi\,dx\nonumber\\
&&=\lim_{\epsilon_k\rightarrow 0} \int_{D}\int_{\mathcal{H}_1(0)}|\xi|J(|\xi|)F_{\epsilon_k|\xi|}'(D_{e}^{{\epsilon_k}|\xi|,-}u^{\epsilon_k})D_{e}^{\epsilon_k|\xi|}\psi\,d\xi\,dx.
\label{functiotolimitminus}
\end{eqnarray}

To finish the proof we identify the limit of the right hand side of \eqref{functiotolimitminus}.
Set $h_{\epsilon_k}=D_{e}^{\epsilon_k |\xi|,-}u^{\epsilon_k}$ and apply H\'older's inequality to find
\begin{eqnarray}
&&\int_{D\times\mathcal{H}_1(0)}|\xi|J(|\xi|)\left(F_{\epsilon_k|\xi|}'(h_{\epsilon_k}) -2f'(0)h_{\epsilon_k}\right )D_e^{\epsilon_k|\xi|}\psi\,d\xi\,dx\nonumber\\
&&\leq C\int_{D\times\mathcal{H}_1(0)}|\xi|J(|\xi|)\left|F_{\epsilon_k|\xi|}'(h_{\epsilon_k}) -2f'(0)h_{\epsilon_k}\right|\,d\xi\,dx\Vert\nabla\psi\Vert_{L^\infty(D)}
\label{firstestimate}
\end{eqnarray}
We estimate the first factor in \eqref{firstestimate} and
apply \eqref{diffeq}, H\"older's inequality, \eqref{L2bound}, and \eqref{usefulbound} to obtain
\begin{eqnarray}
&&\int_{D\times\mathcal{H}_1(0)}|\xi|J(|\xi|)\left |F_{\epsilon_k|\xi|}'(h_{\epsilon_k}) -2f'(0)h_{\epsilon_k}\right |\,d\xi\,dx\nonumber\\
&&\leq\int_{D\times\mathcal{H}_1(0)}|\xi|J(|\xi|)\left |f'(\epsilon_k |\xi||h_{\epsilon_k}|^2) -2f'(0)\right |\left | h_{\epsilon_k}\right |\,d\xi\,dx
\nonumber\\
&&\leq C\left(\nu(D\times\mathcal{H}_1(0))\times |f'(A^2\overline{r}^2)-f'(0)|^2+\epsilon_k\frac{4|f'(0)|^2}{A^2\overline{r}^2}\right).
\label{usefulLemmaA}
\end{eqnarray}
Lemma \ref{twolimitsA} follows on passing to the $\epsilon_k$  zero limit in \eqref{usefulLemmaA} and noting that the choice of $0<A<1$ is arbitrary.

We now prove Lemma \ref{twolimitsB}.
For $\tau>0$ sufficiently small define $K^\tau\subset D$  by  $K^{\tau}=\{x\in D:\,dist(x,S_{u^0(t)})<\tau\}$.  From Hypothesis \ref{remark555} the collection of centroids associated with unstable neighborhoods $C_{\delta,t}$ lie inside $K^\tau$ for $\delta$ sufficiently small. (Otherwise the collection $C_{\delta,t}$ would concentrate about a component of $C_{0,t}$ outside $K^\tau$; contradicting the hypothesis that $S_{u^0(t)}=C_{0,t}$). The collection of all points belonging to unstable neighborhoods associated with centroids in $C_{\delta,t}$  is easily seen to be contained in the slightly larger set $K^{\tau,\delta}=\{x\in \,D; dist(x,K^\tau)<\delta\}$. From Hypothesis \ref{remark556} we may choose test functions $\varphi\in C_0^1(D\setminus K^{\tau,\delta})$ such that for $\epsilon_k$ sufficiently small
\begin{eqnarray}
D_e^{\epsilon_k |\xi|,-}u^{\epsilon_k}\varphi=D_e^{\epsilon_k|\xi|}u^{\epsilon_k}\varphi.
\label{identical}
\end{eqnarray}

We form the test functions $\phi(x,\xi)=\varphi(x)\psi(\xi)$, with $\varphi\in C_0^1(D\setminus K^{\tau,\delta})$ and 
$\psi\in C(\mathcal{H}_1(0))$. From \eqref{L2bound} we may pass to a subsequence  to find that $D_e^{\epsilon_k |\xi|,-}u^{\epsilon_k}$ weakly converges to the limit 
$g(x,\xi)$ in $L^2(D\times\mathcal{H}_1(0);\nu)$. 
With this in mind we write
\begin{eqnarray}
&&\int_{D\times\mathcal{H}_1(0)} g(x,\xi)\phi(x,\xi)|\xi|J(|\xi|)\,d\xi\,dx\nonumber\\
&&=\lim_{\epsilon_k\rightarrow 0}\int_{D\times\mathcal{H}_1(0)}D_e^{\epsilon_k |\xi|,-}u^{\epsilon_k}(x) \phi(x,\xi)|\xi|J(|\xi|)\,d\xi\,dx\nonumber\\
&&=\lim_{\epsilon_k\rightarrow 0}\int_{D\times\mathcal{H}_1(0)}D_e^{\epsilon_k |\xi|}u^{\epsilon_k}(x) \phi(x,\xi)|\xi|J(|\xi|)\,d\xi\,dx\nonumber\\
&&=\lim_{\epsilon_k\rightarrow 0}\int_{D\times\mathcal{H}_1(0)}u^{\epsilon_k}(x)(D_{-e}^{\epsilon_k |\xi|}\varphi(x))\psi(\xi)|\xi|J(|\xi|)\,d\xi\,dx.\label{middlepart}
\end{eqnarray}
Noting that $D_{-e}^{\epsilon_k |\xi|}\varphi(x)$ converges uniformly to $-e\cdot\nabla\varphi(x)$ and from the strong convergence of $u^{\epsilon_k}$ to $u^0$ in $L^2$ we obtain
\begin{eqnarray}
&&=\lim_{\epsilon_k\rightarrow 0}\int_{D\times\mathcal{H}_1(0)}u^{\epsilon_k}(x)(D_{-e}^{\epsilon_k |\xi|}\varphi(x))\psi(\xi)|\xi|J(|\xi|)\,d\xi\,dx\nonumber\\
&&=-\int_{D\times\mathcal{H}_1(0)}u^{0}(x) (e\cdot\nabla\varphi(x))\psi(\xi)|\xi|J(|\xi|)\,d\xi\,dx\nonumber\\
&&=-\int_{D}u^0(x)\,div\left(\varphi(x) \int_{\mathcal{H}_1(0)} e\psi(\xi)|\xi|J(|\xi|)\,d\xi\right)\,dx\nonumber\\
&&=\int_{D}\nabla u^0(x)\cdot\left(\varphi(x) \int_{\mathcal{H}_1(0)} e\psi(\xi)|\xi|J(|\xi|)\,d\xi\right)\,dx\nonumber\\
&&=\int_{D\times{\mathcal{H}_1(0)} }\nabla u^0(x)\cdot e\varphi(x)\psi(\xi)|\xi|J(|\xi|)\,d\xi\,dx,
\label{identifyweakl2}
\end{eqnarray}
where we have made use of  $Du^0\lfloor D\setminus K^{\tau,\delta}=\nabla u^0\,dx$ on the third line of \eqref{identifyweakl2}.
From the density of the span of the test functions we conclude that $g(x,\xi)=\nabla u^0\cdot e$ almost everywhere on $D\setminus K^{\tau,\delta}\times\mathcal{H}_1(0)$.  Since $K^{\tau,\delta}$ can be chosen to have arbitrarily small measure with vanishing $\tau$ and $\delta$ we conclude that  $g(x,\xi)=\nabla u^0\cdot e$ on $D\times\mathcal{H}_1(0)$ a.e. and Lemma \ref{twolimitsB} is proved.

\subsubsection{ Proof of Theorem \ref{bondunstable}}
\label{proofbondunstable}
The set $K^{+,\epsilon_k}$ defined by \eqref{suppset} has the equivalent description given by
\begin{eqnarray}
K^{+,\epsilon_k}=\{(x,\xi)\in D\times\mathcal{H}_1(0);\,|u^{\epsilon_k}(x+\epsilon_k\xi)-u^{\epsilon_k}(x)|>\overline{\eta}\}
\label{equivdescr}
\end{eqnarray}
where $\overline{\eta}$ is the critical stretch given by $\overline{\eta}=\sqrt{\epsilon_k|\xi|}\overline{r}$.
We rewrite the lefthand side of the inequality \eqref{suppsetupper} as
\begin{eqnarray}
\int_{K^{+,\epsilon_k} }|\xi|J(|\xi|)\,d\xi\,dx=\int_{D}\left(\int_{\mathcal{H}_1(0)}\chi^{+,\epsilon_k}(x,\xi)|\xi|J(|\xi|)\,d\xi\right)\,dx,
\label{rewrite}
\end{eqnarray}
where for each $x\in D$, $\chi^{+,\epsilon_k}(x,\xi)$ is defined to be the indicator function given by
\begin{eqnarray}
&&\chi^{+,\epsilon_k}(x,\xi)=1,\hbox{ for }\xi\in\mathcal{H}_1(0);\, |u^{\epsilon_k}(x+\epsilon_k\xi)-u^{\epsilon_k}(x)|>\overline{\eta}\nonumber\\
&&\chi^{+,\epsilon_k}(x,\xi)=0, \hbox{ otherwise. }
\label{indicator}
\end{eqnarray}
Making the change of variable $y=\epsilon_k\xi+x$ the inner integral on the right hand side of \eqref{rewrite} is given by
\begin{eqnarray}
\int_{\mathcal{H}_1(0)}\chi^{+,\epsilon_k}(x,\xi)|\xi|J(|\xi|)\,d\xi=m\times P\left(\{y\in\mathcal{H}_{\epsilon_k}(x);\,|u^{\epsilon_k}(y)-u^{\epsilon_k}(x)|>\overline{\eta}\}\right)
\label{probnewvariable}
\end{eqnarray}
Recall that \eqref{firstb}  shows that the inequality \eqref{suppsetupper} is  uniform both in time and in the length scale of the horizon $\epsilon_k$. This follows from the  uniform bound on the peridynamic potential given by Theorem \ref{Gronwall}.  
Application of \eqref{suppsetupper} gives 
\begin{eqnarray}
\int_{D} P\left(\{y\in\mathcal{H}_{\epsilon_k}(x);\,|u^{\epsilon_k}(t,y)-u^{\epsilon_k}(t,x)|>\overline{\eta}\}\right)\,dx\leq C\epsilon_k.
\label{probnewvariablebound}
\end{eqnarray}
For $A>0$,  Tchebyshev's inequality gives
\begin{eqnarray}
&&A |\{x\in D;\, P\left(\{y\in\mathcal{H}_{\epsilon_k}(x);\,|u^{\epsilon_k}(t,y)-u^{\epsilon_k}(t,x)|>\overline{\eta}\}\right)>A\}|\nonumber\\
&&\leq\int_{D} P\left(\{y\in\mathcal{H}_{\epsilon_k}(x);\,|u^{\epsilon_k}(t,y)-u^{\epsilon_k}(t,x)|>\overline{\eta}\}\right)\,dx.
\label{probnewvariableboundchebyshcv}
\end{eqnarray}
Choosing $A=\sqrt{\epsilon_k}$ and applying \eqref{probnewvariablebound} delivers
\begin{eqnarray}
|\{x\in D;\, P\left(\{y\in\mathcal{H}_{\epsilon_k}(x);\,|u^{\epsilon_k}(t,y)-u^{\epsilon_k}(t,x)|>\overline{\eta}\}\right)>\sqrt{\epsilon_k}\}|<C\sqrt{\epsilon_k}.
\label{ebound}
\end{eqnarray}
Here $C$ is a constant independent of $t$ and $\epsilon_k$.
The collection of centroids $x$ for neighborhoods $\mathcal{H}_{\epsilon_k}(x)$ associated with the instability condition given by 
\begin{eqnarray}
P\left(\{y\in\mathcal{H}_{\epsilon_k}(x);\,|u^{\epsilon_k}(t,y)-u^{\epsilon_k}(t,x)|>\overline{\eta}\}\right)>\sqrt{\epsilon_k}
\label{eboundcond}
\end{eqnarray}
is denoted by $U_{\epsilon_k,t}$. Choose $\epsilon_k=\frac{1}{2^k}$ and \eqref{ebound} imples $|U_{\epsilon_k,t}|<C\frac{1}{\sqrt{2}^k}$.
The unstable set defined by \eqref{unstabelll} is written as
\begin{eqnarray}
C_{\delta,t}=\cup_{\epsilon_k<\delta}U_{\epsilon_k,t}
\label{unstablecollection}
\end{eqnarray}
and from the geometric series we find
\begin{eqnarray}
|C_{\delta,t}|<C\sqrt{\delta}.
\label{boundcdelta}
\end{eqnarray}
and Theorem \ref{bondunstable} follows.

\setcounter{equation}{0} \setcounter{theorem}{0} \setcounter{lemma}{0}\setcounter{proposition}{0}\setcounter{remark}{0}\setcounter{remark}{0}
\setcounter{definition}{0}\setcounter{hypothesis}{0}

\section{Acknowlegements}
\label{Acknowlegements}
The author would like to thank Stewart Silling, Richard Lehoucq and Florin Bobaru for stimulating and fruitful discussions.
This research is supported by NSF grant DMS-1211066, AFOSR grant FA9550-05-0008, and NSF EPSCOR Cooperative Agreement No. EPS-1003897 with additional support from the Louisiana Board of Regents.


\begin{thebibliography}{99}
\bibitem{Ambrosiobook}
L. Ambrosio, N. Fusco, and D. Pallara. {\em Functions of Bounded Variation and Free Discontinuity Problems.} Oxford Mathematical Monographs. Clarendon Press, Oxford, UK, 2000.
\bibitem{DeGiorgiAmbrosio}
L. Ambrosio and E. De Giorgi. {\em Un nuovo tipo di funzionale del calcolo delle variazioni}, Atti della Accademia Nazionale dei Lincei, Rendiconti della Classe di Scienze Fisiche, Matematiche e Naturali.  {\bf 82} (1989) 199--210.
\bibitem{Ambrosio}
L. Ambrosio. {\em A compactness theorem for a special class of functions of bounded variation.} Boll. Un. Mat. It. {\bf 3-B} (1989) 857--881.
\bibitem{AmbrosioBrades}
L. Ambrosio and A. Brades. {\em Energies in SBV and variational models in fracture mechanics.} Homogenization and Applications to Materials Science {\bf 9}. D. Cioranescu, A. Damlamian, P. Donnato eds, Gakuto, Gakkotosho, Tokyo, Japan (1997)  1--22.

\bibitem{Bobaru1}
F. Bobaru and W. Hu. {\em The meaning, selection, and use of the Peridynamic horizon and its relation to crack branching in brittle materials.} International Journal of Fracture {\bf 176} (2012) 215-Ð222.
\bibitem{Hughes}
M. Borden, C. Verhoosel, M. Scott, T. Hughes, and C. Landis. {\em A phase-field description of dynamic brittle fracture.} Computer Methods in Applied Mechanics and Engineering {\bf 217-220} (2012) 77-95.
\bibitem{BourdinFrancfortMarigo}
B. Bourdin, G. Francfort, and J.-J. Marigo. {\em The variational approach to fracture}. J. Elasticity {\bf 91} (2008) 5--148.
\bibitem{BourdinLarsenRichardson}
B. Bourdin, C. Larsen, C. Richardson. {\em A time-discrete model for dynamic fracture based on crack regularization.} Int. J. Fract  {\bf 168}  (2011) 133--143.
\bibitem{Braides}
A. Braides. {\em Discrete approximation of functionals with jumps and creases.} Homogenization, 2001 (Naples), Gakuto Internat. Ser. Math. Sci. Appl. {\bf 18} Gakkotosho, Tokyo, 2003, 147--153.
\bibitem{Buttazzo}
H. Attouch, G. Buttazzo, and G. Michaille. {\em Variational Analysis in Sobolev and BV Spaces: applications to PDEs and optimization.} MPS-SIAM series on optimization. SIAM, Philadelphia, PA, 2006.
\bibitem{LarsenDalMaso}
G. Dal Maso and C. J. Larsen. {\em  Existence for wave equations on domains with
arbitrary growing cracks.} Rend. Lincei Mat. Appl. {\bf 22} (2011) 387--408.
\bibitem{Driver}
B. Driver. {\em Analysis Tools With Applications} E-book, Springer, Berlin, 2003.
\bibitem{BhattacharyaDyal}
K. Dyal and K. Bhattacharya. {\em Kinetics of phase transformations in the perydanimic formulation of continuum mechanics.} J. Mech. Phys. Solids. {\bf 54} (2006) 1811-1842.
\bibitem{DuGunzbergerlehoucqmengesha}
Q. Du, M. Gunzburger, R. Lehoucq, and K. Zhou. {\em A nonlocal vector calculus, nonlocal volume-constrained problems, and nonlocal balance laws.} (Mathematical Models and Methods in Applied Sciences (M3AS), {\bf 23} (2013) 493--540.
\bibitem{EmmrichWeckner}
E. Emmrich and O. Weckner.  {\em On the well-posedness of the linear peridynamic model and its convergence towards the Navier equation of linear elasticity.} Communications in Mathematical Sciences {\bf 5} (2007) 851--864.
\bibitem{FrancfortMarigo}
G. Francfort and J.-J. Marigo. {\em Revisiting brittle fracture as an energy minimization problem.} J. Mech. Phys. Solids {\bf 46} (1998) 1319--1342.
\bibitem{Freund}
L.B. Freund. {\em Dynamic Fracture Mechanics.} Cambridge Monographs on Mechanics and  Applied Mathematics. Cambridge University Press, Cambridge, UK, 1998.
\bibitem{Evans}
L. C. Evans. {\em Partial Differential Equations.} Graduate Studies in Mathematics Vol. 19, American Mathematical Society, Providence, RI, 2010.
\bibitem{SillingAscari3}
W. Gerstle, N. Sau, and S. Silling. {\em Peridynamic Modeling of Concrete Structures.} Nuclear Engineering and Design {\bf 237} (2007) 1250--1258.
\bibitem{Gobbino1}
M. Gobbino. {\em Finite difference approximation of the Mumford-Shah Functional.} Comm. Pure appl. Math. {\bf 51} (1998) 197--228.
\bibitem{Gobbino2}
M. Gobbino. {\em Gradient flow for the one-dimensional Mumford-Shah functional.} Ann. Scuola Norm. Sup. Pisa Cl. Sci. (4) Vol. XXVII (1998) 145--193.
\bibitem{Gobbino3}
M. Gobbino and M.G. Mora. {\em Finite difference approximation of free discontinuity problems.} Royal Soc. Edinburgh Proceedings A. {\bf 131} (2001) 567--595.
\bibitem{Bobaru2}
W. Hu, YD. Ha, and F. Bobaru. {\em Modeling dynamic fracture and damage in a fiber-reinforced composite lamina with peridynamics.} International Journal for Multiscale Computational Engineering {\bf 9} (2011) 707--726.
\bibitem{Larsen}
C. J. Larsen. {\em Models for dynamic fracture based on Griffith's criterion,}
in IUTAM Sympo-sium on Variational Concepts with Applications to the Mechanics of
Materials (Klaus Hackl, ed.), Springer, 2010, pp. 131--140.
\bibitem{LarsenOrtiner}
C. Larsen, C. Ortner, and E. Suli. {\em Existence of solutions to a regularized model of dynamic fracture.} M3AS {\bf 20} (2010) 1021--1048.
\bibitem{MumfordShah}
D. Mumford and J. Shah. {\em Optimal approximation by piecewise smooth functions and associated variational problems.} Comm. Pure Appl. Math. {\bf 17} (1989) 577--685.
\bibitem{Ortiz}
B. Schmidt, F. Fraternali, and M. Ortiz. {\em Eigenfracture: an eigendeformation approach to variational fracture.} Multiscale Model. Simul. {\bf 7} (2009) 1237--1266.
\bibitem{Silling1}
S. A. Silling. {\em Reformulation of Elasticity Theory for Discontinuities and Long-Range Forces.} J. Mech. Phys. Solids {\bf 48} (2000) 175Ð209.
\bibitem{SillingAscari2}
S. A. Silling and E. Askari. {\em A meshfree method based on the peridynamic model of solid mechanics.} Computers and Structures {\bf 83} (2005) 1526--1535.
\bibitem{SillingBobaru}
S. A. Silling and F. Bobaru. {\em Peridynamic Modeling of Membranes and Fibers.} International Journal of Non-Linear Mechanics {\bf 40} (2005) 395--409.
\bibitem{SillingWecknerBobaru}
S. Silling, O. Weckner, E. Askari, and F. Bobaru. {\em Crack nucleation in a peridynamic solid.} International Journal of Fracture {\bf 162} (2010) 219--227.
\bibitem{SillingLehoucq}
S. Silling and R. Lehoucq. {\em Convergence of Peridynamics to Classical Elasticity Theory.} Journal of Elasticity  {\bf 93} (2008) 13--37.
\bibitem{WecknerAbeyaratne}
O. Weckner and R. Abeyaratne. {\em The Effect of Long-Range Forces on the Dynamics of a Bar.} Journal of the Mechanics and Physics of Solids {\bf 53} (2005) 705--728.
\end{thebibliography}
\end{document}